\documentclass[a4paper,12pt]{article} 

\usepackage{color}

\usepackage{amsfonts, amsmath, amsthm, amssymb}
\usepackage[T1]{fontenc}
\usepackage[cp1250]{inputenc}
\usepackage{xcolor}
\usepackage{graphicx}
\usepackage{amssymb}
\usepackage{amsmath}
\usepackage{mathptmx}
\usepackage{helvet}
\usepackage{courier}
\usepackage{txfonts}
\usepackage{tikz} 
\usetikzlibrary{arrows}
\usepackage{type1cm}

\usepackage{verbatim}

\usepackage{graphicx}
\usepackage{epsfig,amscd,amssymb,amsxtra,amsmath,amsthm}
\usepackage{type1cm}
\usepackage[T1]{fontenc}
\usepackage{graphics}
\usepackage[mathscr]{eucal}
\usepackage[all]{xy}
\usepackage{amsmath,amscd}

\newcommand{\orbit}{\mathcal{O}^{\oplus}}

\newcommand{\Orbit}{\mathcal{O}_G^{\oplus}}
\newcommand{\orbitt}{\mathcal{U}^{\oplus}}
\newcommand{\tr}{\textup{tr}}
\newcommand{\intr}{\textup{intr}}
\newcommand{\illegal}{\textup{illegal}}
\newcommand{\legal}{\textup{legal}}
\newcommand{\isolated}{\textup{isolated}}
\newcommand{\level}{\textup{level}}
\newcommand{\tplus}{T^{+}}

\newcommand{\trans}{\textup{trans}}
\newcommand{\intrans}{\textup{intrans}}

\newtheorem{theorem}{Theorem}[section]

\newtheorem{definition}[theorem]{Definition}
\newtheorem{lemma}[theorem]{Lemma}

\newtheorem{example}[theorem]{Example}

\newtheorem{corollary}[theorem]{Corollary}

\newtheorem{observation}[theorem]{Observation}

\newtheorem*{nonumclaim}{Claim}

\newcommand{\Cl}  {\mathop{\rm Cl}\nolimits}

\begin{document}

\def\joinrel{\mkern-3mu}
\newcommand{\varproj}{\displaystyle \lim_{\multimapinv\joinrel-\joinrel-}}

\title{Transitive  points in CR-dynamical systems}
\author{Iztok Bani\v c, Goran Erceg, Sina Greenwood, Judy Kennedy}
\date{}

\maketitle

\begin{abstract}
\noindent We study different types of transitive points in CR-dynamical systems $(X,G)$ with closed relations   $G$ on compact metric spaces $X$.  We also introduce transitive and dense orbit transitive CR-dynamical systems  and discuss  their properties and the relations between them.  This generalizes the notion of transitive topological dynamical systems $(X,f)$.
\end{abstract}
\-


\section{Introduction}\label{s000}

In topological dynamical systems theory,  the study of the chaotic behaviour of a dynamical system is  often based on some  properties of topological spaces or some properties of continuous functions.  One of the commonly studied properties in the theory of topological dynamical systems is the transitivity of a dynamical system $(X,f)$ or the transitivity of the function $f$.  This often reduces to studying transitive or intransitive points of dynamical systems.  Regarding this, the study of forward orbits and backward orbits of points is also required.  Backward orbits of points are actually forward orbits of points in the dynamical system $(X,f^{-1})$, if $f^{-1}$ is well-defined.  But usually, $f^{-1}$ is not a well-defined function, therefore, a more general tool is needed to study these properties. Note that for a continuous function $f:X\rightarrow X$,  the set
$$
\Gamma(f)^{-1}=\{(y,x) \in X\times X \ | \  y=f(x)\}
$$  
is a closed relation on $X$ that describes best the dynamics of $(X,f)$ in the backward direction  when $f^{-1}$ is not well-defined.  So,  generalizing topological dynamical systems $(X,f)$ to topological dynamical systems $(X,G)$ with closed relations $G$ on $X$ by making the identification $(X,f)=(X,\Gamma(f))$ is only  natural.

In this paper, we generalize the notion of transitivity from topological dynamical systems to topological dynamical systems with closed relations.  We introduce notions of different types of transitive points, different types of dense orbit transitivity as well as the transitivity of dynamical systems with closed relations.   In the past,  many similar generalizations of topological or dynamical objects have already been introduced. One of them  was presented in 2004 by Ingram and Mahavier  \cite{ingram,mah} introducing inverse limits of inverse sequences of compact metric spaces $X$ with upper semi-continuous set-valued bonding functions $f$. Their graphs $\Gamma(f)$ are examples of closed relations on $X$.  These inverse limits provide a valuable extension to the role of inverse limits in the study of dynamical systems and continuum theory.  For example,   Kennedy and Nall  have developed a simple method for constructing families of $\lambda$-dendroids \cite{KN}. Their method involves inverse limits of inverse sequences with upper semi-continuous set-valued functions on closed intervals with simple bonding functions. 
Such generalizations have proven to be  useful (also in applied areas); frequently,  when  constructing a model for  empirical data, continuous (single-valued) functions fall short, and the data are better modelled by upper semi-continuous set-valued functions, or sometimes, even closed relations that are not set-valued functions are required. The Christiano-Harrison model   from macroeconomics is one such example \cite{christiano}.  The study of inverse limits of inverse sequences with upper semi-continuous set-valued  functions is rapidly gaining momentum - the recent books by Ingram \cite{ingram-knjiga}, and by Ingram and Mahavier \cite{ingram}, give a comprehensive exposition of this research prior to 2012.   

Also,  several papers on the topic of dynamical systems with (upper semi-continuous) set-valued functions have appeared recently, see \cite{BEGK,CP,LP,LYY,LWZ,KN,KW,MRT,R,SS,SS2}, where more references may be found.  However,  there is not much known of such dynamical systems and therefore,  there are many properties of such set-valued dynamical systems that are yet to be studied.  We proceed as follows. In sections that follow Section \ref{s00}, where basic definitions are given, we discuss the following topics:
\begin{itemize}
\item Transitive and intransitive points in CR-dynamical systems (Section \ref{s1}).
\item 3-transitive points that are not 2-transitive (Section \ref{s2}).
\item Transitivity trees (Section \ref{s3}).
\item Dense orbit transitive CR-dynamical systems  (Section \ref{s4}).
\item Transitive CR-dynamical systems (Section \ref{s5}).
\end{itemize}
In each section when introducing a new object in CR-dynamical systems,  we first revisit transitive  dynamical systems $(X,f)$ and then, we  generalize this property from dynamical systems $(X,f)$ to dynamical systems with closed relations $(X,G)$ by making the identification $(X,f)=(X,\Gamma(f))$.   Results about dynamical systems $(X,f)$,  presented in this paper are well-known.  Therefore, we omit them. The reader can track their proofs by a little help from S.~Kolyada's and L.~Snoha's paper ``Topological Transitivity'' \cite{KS}, where a wonderful overview of transitive dynamical systems is given,  or by E.~Akin's book ``General Topology of Dynamical Systems'' \cite{A}, where dynamical systems using closed relations are presented. 

At the end of the paper, two examples are given proving that the statement of \cite[Theorem 9, page 3]{SS},  saying that any dynamical system with an upper semi-continuous set-valued function is transitive if and only if there is a point with a dense orbit,  is incorrect.

\section{Definitions and notation} \label{s00}

In this section, basic definitions and well-known results that are needed later in the paper are presented. All spaces in this paper are non-empty compact metric spaces.
\begin{definition}
Let $X$ and $Y$ be metric spaces, and let $f:X\rightarrow Y$ be a function.  We use  
$$
\Gamma(f)=\{(x,y)\in X\times Y \ | \ y=f(x)\}
$$
to denote \emph{ \color{blue} the graph of the function $f$}.
\end{definition}
\begin{definition}
If $X$ is a compact metric space, then \emph{ \color{blue}  $2^X$ }denotes the set of all  non-empty closed subsets of $X$.  
\end{definition}
%

\begin{definition}
Let $X$ be a compact metric space and let $G\subseteq X\times X$ be a relation on $X$. If $G\in 2^{X\times X}$, then we say that $G$ is \emph{ \color{blue} a closed relation on $X$}.  
\end{definition}

\begin{definition}
Let $X$  be a set and let $G$ be a relation on $X$.  Then we define  
$$
G^{-1}=\{(y,x)\in X\times X \ | \ (x,y)\in G\}
$$
to be \emph{ \color{blue} the inverse relation of the relation $G$} on $X$.
\end{definition}
\begin{definition}
Let $X$ be a compact metric space and let $G$ be a closed relation on $X$. Then we call
$$
\star_{i=1}^{m}G=\Big\{(x_1,x_2,x_3,\ldots ,x_{m+1})\in \prod_{i=1}^{m+1}X \ | \ \textup{ for each } i\in \{1,2,3,\ldots ,m\}, (x_{i},x_{i+1})\in G\Big\}
$$
for each positive integer $m$, \emph{ \color{blue} the $m$-th Mahavier product of $G$}, and
$$
\star_{i=1}^{\infty}G=\Big\{(x_1,x_2,x_3,\ldots )\in \prod_{i=1}^{\infty}X \ | \ \textup{ for each positive integer } i, (x_{i},x_{i+1})\in G\Big\}
$$
\emph{ \color{blue} the infinite  Mahavier product of $G$}.
\end{definition}
\begin{observation}
Let $X$ be a compact metric space, let $f:X\rightarrow X$ be a continuous function. 
Then 
$$
\star_{n=1}^{\infty}\Gamma(f)^{-1}=\varprojlim(X,f).
$$
\end{observation}

%
%
%
\begin{definition}
Let $X$ be a compact metric space.  For each positive integer $k$,  we use $\pi_k:\prod_{i=1}^{\infty}X\rightarrow X$ to denote \emph{ \color{blue}  the $k$-th standard projection} from $\prod_{i=1}^{\infty}X$ to $X$. We also use use $\pi_k:\prod_{i=1}^{n}X\rightarrow X$ to denote \emph{ \color{blue}  the $k$-th standard projection} from $\prod_{i=1}^{n}X$ to $X$, where $n$ is a positive integer and $k\in \{1,2,3,\ldots, n\}$. 
\end{definition}
\begin{definition}
Let $G$ be a relation on $X$, let $A\subseteq X$ and let $x\in X$.  Then we define 
$$
G(x)=\{y\in X \ | \ (x,y)\in G\}
$$
and 
$$
G(A)=\bigcup_{x\in A}G(x).
$$ 
For each positive integer $n$, we also define
$$
G^n(x)=\{y\in X \ | \ \textup{ there is } \mathbf x\in \star_{i=1}^{n}G \textup{ such that } \pi_1(\mathbf x)=x \textup{ and } \pi_{n+1}(\mathbf x)=y\}
$$
and 
$$
G^n(A)=\bigcup_{x\in A}G^n(x).
$$
For each positive integer $n$, we also define
$$
G^{-n}(x)=(G^{-1})^n(x)
$$
and 
$$
G^{-n}(A)=\bigcup_{x\in A}G^{-n}(x).
$$
Finally, we define
$$
G^{\omega}(X)=\bigcap_{n=1}^{\infty}G^{n}(X)
$$
and
$$
G^{-\omega}(X)=\bigcap_{n=1}^{\infty}G^{-n}(X).
$$
\end{definition}
We conclude this section by stating the following observations.
\begin{observation}\label{encica}
Let $(X,G)$ be a CR-dynamical system,  let $x\in X$ and let $n$ be a positive integer. Then $G^{n+1}(x)=G(G^n(x))$.
\end{observation}
\begin{observation}\label{dvojcica}
Let $(X,G)$ be a CR-dynamical system,  let $x,y\in X$, and let $n$ be a non-negative integer.  The following statements are equivalent.
\begin{enumerate}
\item $x\in G^n(y)$. 
\item There is $\mathbf x=(x_1,x_2,x_3,\ldots, x_{n},x_{n+1})\in \star_{i=1}^{n}G$ such that 
$$
x_1=y \textup{ and } x_{n+1}=x.
$$
\end{enumerate} 
\end{observation}
\begin{observation}\label{trojcica}
Let $(X,G)$ be a CR-dynamical system,  let $x,y\in X$, and let $n$ be a non-negative integer.  The following statements are equivalent.
\begin{enumerate}
\item $x\in G^{-n}(y)$. 
\item There is $\mathbf x=(x_1,x_2,x_3,\ldots, x_{n},x_{n+1})\in \star_{i=1}^{n}G^{-1}$ such that 
$$
x_1=y \textup{ and } x_{n+1}=x.
$$
\item There is $\mathbf x=(x_1,x_2,x_3,\ldots, x_{n},x_{n+1})\in \star_{i=1}^{n}G$ such that 
$$
x_1=x \textup{ and } x_{n+1}=y.
$$
\end{enumerate} 
\end{observation}
\section{Transitive and intransitive points in CR-dynamical systems}\label{s1}
First, we revisit dynamical systems, and then,  transitive  and intransitive points in dynamical systems. 
\begin{definition}
Let $X$ be a non-empty compact metric space and $f:X\rightarrow X$ a function. If $f$ is continuous, then we say that $(X,f)$ is a dynamical system.
\end{definition}
\begin{definition}
Let $(X,f)$ be a dynamical system and let $x\in X$. The sequence 
$$
\mathbf x=(x,f(x),f^2(x),f^3(x),\ldots)\in \star_{i=1}^{\infty}\Gamma(f)
$$   
is called \emph{ \color{blue} the trajectory of $x$ in $(X,f)$.} The set 
$$
\mathcal O_f^{\oplus}(x)=\{x,f(x),f^2(x),f^3(x),\ldots\}
$$   
is called \emph{ \color{blue} the forward orbit of $x$ in $(X,f)$}.
\end{definition}

\begin{definition}
Let $(X,f)$ be a dynamical system and let $x\in X$.  If $\Cl(\orbit(x))=X$, then $x$ is called \emph{ \color{blue} a transitive point in $(X,f)$}.  Otherwise it is \emph{ \color{blue} an intransitive point in $(X,f)$}. We use \emph{ \color{blue} $\tr(f)$} to denote the set
$$
\tr(f)=\{x\in X \ | \ x \textup{ is a transitive point in } (X,f)\}
$$
and \emph{ \color{blue} $\intr(f)$} to denote the set
$$
\intr(f)=\{x\in X \ | \ x \textup{ is an intransitive point in } (X,f)\}.
$$
\end{definition}
Next, we generalize these to dynamical systems with closed relations. 
\begin{definition}
Let $X$ be a non-empty compact metric space and let $G$ be a closed relation on $X$. We say that $(X,G)$ is \emph{ \color{blue} a dynamical system with a closed relation} or, briefly, \emph{ \color{blue} a CR-dynamical system}.
\end{definition}
\begin{observation}
Let $X$ be a non-empty compact metric space and $f:X\rightarrow X$ a function.  Then  $(X,f)$ is a dynamical system  if and only if $(X,\Gamma(f))$ is a CR-dynamical system.
\end{observation}

\begin{definition}
Let $(X,G)$ be a CR-dynamical system and let $x_0\in X$. We use \emph{ \color{blue}  $T_G^{+}(x_0)$} to denote the set 
$$
T_G^{+}(x_0)=\{\mathbf x\in \star_{i=1}^{\infty}G \ | \  \pi_1(\mathbf x)=x_0\}\subseteq \star_{i=1}^{\infty}G.
$$
\end{definition}
\begin{observation}
Let $(X,f)$ be a dynamical system and let $x_0\in X$.  Then
$$
T_{\Gamma(f)}^{+}(x_0) =\{\mathbf x\},
$$
where $\mathbf x$ is the trajectory of $x_0$ in $(X,f)$.
\end{observation}
\begin{definition}
Let $(X,G)$ be a CR-dynamical system and let $x\in X$. We say that
\begin{enumerate}
\item $x$ is \emph{\color{blue} an illegal point in $(X,G)$}, if $T_G^{+}(x)=\emptyset$. We use 
\emph{\color{blue}$\illegal(G)$} to denote the set of illegal points in $(X,G)$.
\item $x$ is \emph{\color{blue} a legal point in $(X,G)$}, if $T_G^{+}(x)\neq \emptyset$. We use 
\emph{\color{blue}$\legal(G)$} to denote the set of legal points in $(X,G)$.
\end{enumerate}
\end{definition}
\begin{observation}
Let $(X,f)$ be a dynamical system. Then $\illegal(\Gamma(f))=\emptyset$ and $\legal(\Gamma(f))=X$.
\end{observation}
\begin{example}
Let $X=\{1,2\}$ and let $G=\{(1,1)\}$. Then $\illegal(G)=\{2\}$ and $\legal(G)=\{1\}$.
\end{example}
\begin{theorem}\label{main1}
Let $(X,G)$ be a CR-dynamical system and let $x\in X$. The following statements are equivalent.
\begin{enumerate}
\item\label{b1} $x\in \illegal(G)$.
\item\label{b2} There is a positive integer $n$ such that $G^n(x)=\emptyset$.
\item\label{b3} $x\in X\setminus G^{-\omega}(X)$.
\end{enumerate}
\end{theorem}
\begin{proof}
To show that \ref{b1} implies \ref{b2}, let $x\in \illegal(G)$ and suppose that for each positive integer $n$,  $G^n(x)\neq \emptyset$.  For each positive integer $n$, let 
$x_n\in G^n(x)$ and let 
$$
\mathbf x_n\in \star_{i=1}^{n}G
$$
be such that $\pi_1(\mathbf x_n)=x$ and $\pi_{n+1}(\mathbf x_n)=x_n$ (by Observation \ref{dvojcica}, such a point $\mathbf x_n$ does exist). Next,  for each positive integer $n$, let
$$
\mathbf z_n\in  \prod_{i=1}^{\infty}X
$$
be such a point that for each $k\in \{1,2,3,\ldots, n+1\}$, $\pi_k(\mathbf z_n)=\pi_k(\mathbf x_n)$. Let $\mathbf z=(z_1,z_2,z_3,\ldots)\in \prod_{i=1}^{\infty}X$ and let $(\mathbf z_{i_n})$ be a convergent subsequence of the sequence $(\mathbf z_n)$ such that 
$$
\lim_{n\to \infty}\mathbf z_{i_n}=\mathbf z.
$$
Then $\mathbf z\in \star_{i=1}^{\infty}G$ and $\pi_1(\mathbf z)=x$. It follows that $x\in \legal(G)$  -- a  contradiction.

To prove the implication from \ref{b2} to \ref{b3}, let $n$ be a positive integer such that $G^n(x)=\emptyset$ and suppose that $x\in G^{-\omega}(X)$.  Then $x\in G^{-n}(X)$ for each positive integer $n$.  For each positive integer $n$, let 
$$
\mathbf x_n\in \star_{i=1}^{n}G^{-1}
$$
be such a point that $\pi_{n+1}(\mathbf x_n)=x$ (by Observation \ref{trojcica}, such a point $\mathbf x_n$ exists). Also, for each positive integer $n$, let 
$$
\mathbf y_n=(\pi_{n+1}(\mathbf x_n), \pi_{n}(\mathbf x_n), \pi_{n-1}(\mathbf x_n),\ldots ,\pi_{1}(\mathbf x_n) ).
$$
Then  $\mathbf y_n\in \star_{i=1}^{n}G$ for each positive integer $n$. It follows that  for each positive integer $n$, $G^n(x)\neq\emptyset$  -- a  contradiction.

To show the implication from \ref{b3} to \ref{b1}, let $x\in X\setminus G^{-\omega}(X)$ and show that $x\in \illegal(G)$.  Suppose that $x\in \legal(G)$.  Let 
$$
\mathbf x=(x_1,x_2,x_3,\ldots)\in \star_{i=1}^{\infty}G
$$
be such that $x_1=x$.   For each positive integer $n$, let 
$$
\mathbf x_n=(x_{n+1},x_n,x_{n-1},\ldots, x_3,x_2,x_1).
$$
Then $\mathbf x_n\in \star_{i=1}^{n}G^{-1}$ for each positive integer $n$.  Therefore, for each positive integer $n$, $x\in G^{-n}(X)$. It follows that $x\in G^{-\omega}(X)$  -- a  contradiction. 
\end{proof}
%
\begin{lemma}\label{closed1}
Let $(X,G)$ be a dynamical system  and let $A$ be a closed subset of $X$. Then for each positive integer $n$, $G^n(A)$ and $G^{-n}(A)$ are closed in $X$. 
\end{lemma}
\begin{proof}
To prove that for each positive integer $n$, $G^n(A)$ is closed in $X$, we use induction on $n$.  First, we prove that $G(A)$ is closed in $X$.  Let $p_2:X\times X\rightarrow X$ be the second standard projection, $p_2(x,y)=y$ for any $(x,y)\in X\times X$.  Note that 
$$
G(A)=p_2(G\cap (A\times X)).
$$ 
Since $G\cap (A\times X)$ is compact and $p_2$ is continuous, it follows that $p_2(G\cap (A\times X))$ is compact. Therefore, $G(A)$ is closed in $X$.

Next, let $n$ be a positive integer and suppose that $G^n(A)$ is closed in $X$.  By the previous step of the proof,  it follows that $G(G^n(A))$ is closed in $X$.  To finish the proof, note that it follows from Observation \ref{encica} that $G(G^n(A))=G^{n+1}(A)$.
We have just proved that for each positive integer $n$, $G^n(A)$ is closed in $X$. 

Next, let $n$ be a positive integer.  Since $G^{-n}(A)=(G^{-1})^n(A)$ and since also $G^{-1}$ is a closed relation on $X$, it follows that $G^{-n}(A)$ is closed in $X$.
\end{proof}
\begin{corollary}
Let $(X,G)$ be a CR-dynamical system.  Then for each positive integer $n$, 
$$
G^{n+1}(X)\subseteq G^n(X),
$$
so $G^{\omega}(X)$ is a non-empty compactum.
\end{corollary}\begin{proof}
The corollary follows directly from Lemma \ref{closed1}.
\end{proof}
\begin{theorem}\label{illegal_open}
Let $(X,G)$ be a CR-dynamical system.  The set $\illegal(G)$ is open in $X$.
\end{theorem}
\begin{proof}
By Theorem \ref{main1}, 
$$
\illegal(G)=X\setminus G^{-\omega}(X).
$$
It follows from
$$
X\setminus G^{-\omega}(X)=X\setminus \Big(\bigcap_{n=1}^{\infty} G^{-n}(X)\Big)=\bigcup_{n=1}^{\infty}\Big(X\setminus G^{-n}(X)\Big)
$$
that 
$$
\illegal(G)=\bigcup_{n=1}^{\infty}\Big(X\setminus G^{-n}(X)\Big).
$$
It follows from  Lemma \ref{closed1} that for each positive integer $n$,  $G^{-n}(X)$ is closed in $X$,  and  $\illegal(G)$ is open in $X$.
\end{proof}
\begin{corollary}\label{cor}
Let $(X,G)$ be a CR-dynamical system.  The set $\legal(G)$ is closed in $X$.
\end{corollary}
\begin{proof}
Since $\legal(G)=X\setminus \illegal(G)$,  the corollary follows  directly from Theorem \ref{illegal_open}. 
\end{proof}
\begin{definition}
Let $(X,G)$ be a CR-dynamical system, let $\mathbf x\in \star_{i=1}^{\infty}G$, and let $x_0\in X$. We 
 \begin{enumerate}
 \item say that $\mathbf x$ is \emph{\color{blue} a trajectory of $x_0$ in $(X,G)$}, if $\pi_1(\mathbf x)=x_0$. 
 \item use \emph{ \color{blue}  $\mathcal O_G^{\oplus}(\mathbf x)$} to denote the set
 $$
 \mathcal O_G^{\oplus}(\mathbf x)=\{\pi_k(\mathbf x) \ | \  k \textup{ is a positive integer}\}.
 $$
 We call the set $\mathcal O_G^{\oplus}(\mathbf x)$ \emph{\color{blue} a forward orbit of $x_0$ in $(X,G)$}
 \item use \emph{ \color{blue}  $\mathcal U_G^{\oplus}(x_0)$} to denote the set
 $$
 \mathcal U_G^{\oplus}(x_0)=\bigcup_{\mathbf y\in T_G^{+}(x_0)}\mathcal O_G^{\oplus}(\mathbf y).
 $$
 \end{enumerate}
\end{definition}
\begin{definition}
Let $(X,G)$ be a CR-dynamical system and let $x\in \legal(G)$.  We say that 
\begin{enumerate}
\item $x$ is \emph{\color{blue} a type 1 transitive point} or \emph{\color{blue} 1-transitive point} in $(X,G)$, if for each $\mathbf x\in T_G^{+}(x)$, $\orbit_G(\mathbf x)$ is dense in $X$. We use \emph{\color{blue} $\trans_1(G)$ } to denote the set of 1-transitive points in $(X,G)$.
\item $x$ is \emph{\color{blue} a type 2 transitive point} or \emph{\color{blue} 2-transitive point} in $(X,G)$, if there is $\mathbf x\in T_G^{+}(x)$ such that $\orbit_G(\mathbf x)$ is dense in $X$. We use \emph{\color{blue} $\trans_2(G)$ } to denote the set of 2-transitive points in $(X,G)$.
\item $x$ is \emph{\color{blue} a type 3 transitive point} or \emph{\color{blue} 3-transitive point} in $(X,G)$, if $\mathcal U_G^{\oplus}(x)$ is dense in $X$. We use \emph{\color{blue} $\trans_3(G)$ } to denote the set of 3-transitive points in $(X,G)$.
\item $x$ is \emph{\color{blue} an intransitive point} in $(X,G)$, if $x$ is not 3-transitive. We use \emph{\color{blue} $\intrans(G)$ } to denote the set of intransitive points in $(X,G)$.
\end{enumerate}
\end{definition}

\begin{observation}
Let $(X,G)$ be a CR-dynamical system. Note that 
$$
\trans_1(G)\subseteq \trans_2(G)\subseteq \trans_3(G).
$$
\end{observation}
\begin{observation}
Let $(X,G)$ be a CR-dynamical system. Note that 
$$
\legal(G)=\trans_3(G)\cup \intrans(G).
$$
Therefore, by Corollary \ref{cor}, the set $\trans_3(G)\cup \intrans(G)$ is closed in $X$.
\end{observation}
\begin{observation}
Let $(X,f)$ be a dynamical system. Then
$$
\tr(f)=\trans_1(\Gamma(f))=\trans_2(\Gamma(f))= \trans_3(\Gamma(f)).
$$
\end{observation}
In the following example, we show that there are CR-dynamical systems $(X,G)$ that have type 2 transitive points which are not type 1 transitive. 

\begin{example}\label{ex1}
Let $X=[0,1]$ and let 
$$
G=\Big([0,1]\times \Big\{\frac{1}{2}\Big\}\Big)\cup\Big(\Big\{\frac{1}{2}\Big\}\times [0,1]\Big)
$$
see 
Figure \ref{figure1}.
\begin{figure}[h!]
	\centering
		\includegraphics[width=20em]{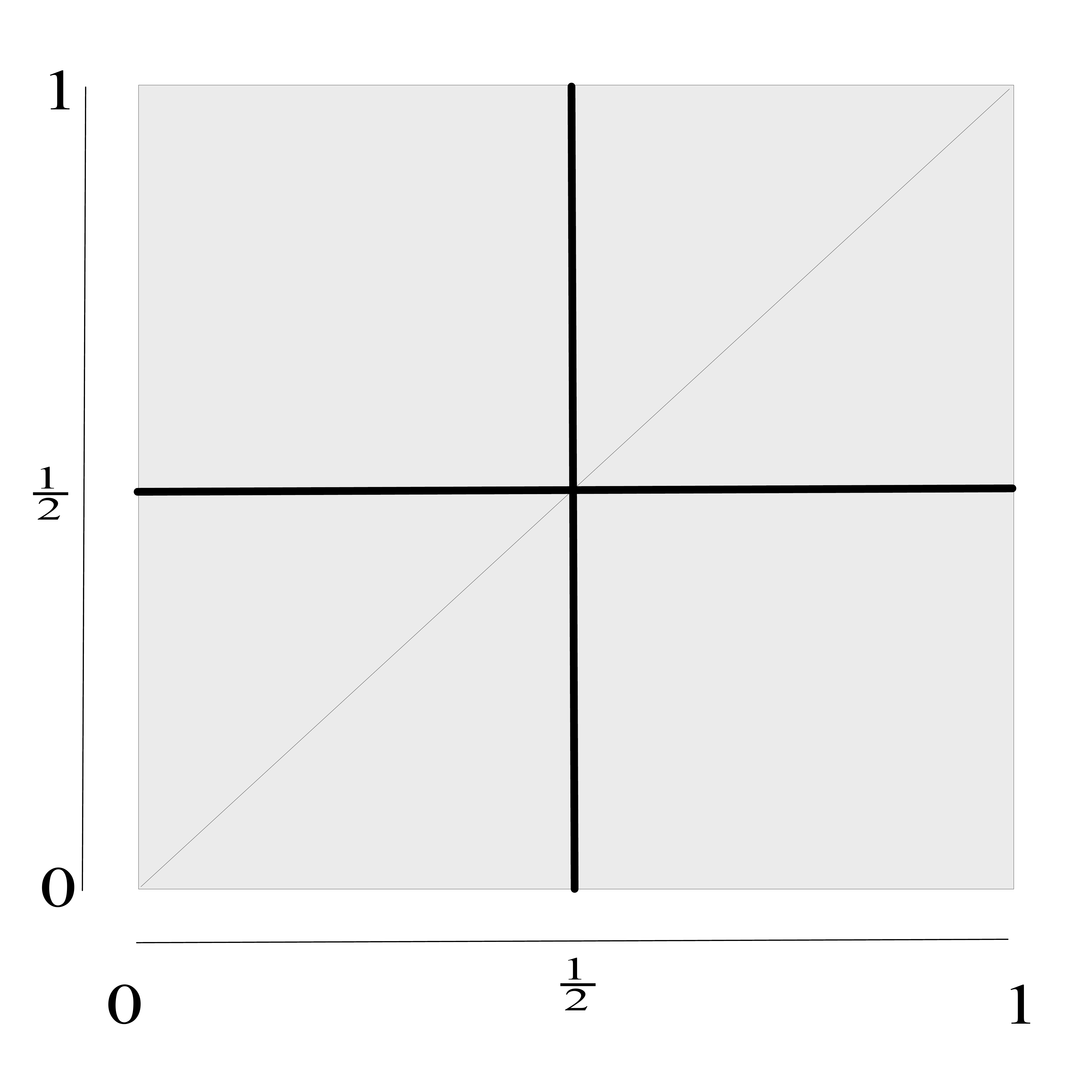}
	\caption{The relation  $G$ from Example \ref{ex1}}
	\label{figure1}
\end{figure} 
\noindent Then $(X,G)$ is a CR-dynamical system. Let $x=\frac{1}{2}$. Then $\mathbf x=(x,x,x,\ldots)\in \tplus_G(x)$ is such a point that $\orbit_G(\mathbf x)=\{x\}$ is not dense in $X$.  Let 
$$
[0,1]\cap \mathbb Q=\{q_1,q_2,q_3,\ldots\}
$$
 be the set of rationals in $[0,1]$, let $x_1=x$,  for each positive integer $n$,  let $x_{2n}=\frac{1}{2}$ and $x_{2n+1}=q_n$, and let $\mathbf x=(x_1,x_2,x_3,\ldots)$. Then $\mathbf x\in T_G^{+}(x)$ is such a point  that $\Cl\Big(\Orbit(\mathbf x)\Big)=X$. Therefore, 
$$
x\in \trans_2(G)\setminus \trans_1(G).
$$
\end{example}

In the following example, we show that there are CR-dynamical systems $(X,G)$ that have type 3 transitive points which are not type 2 transitive. 

\begin{example}\label{ex2}
Let $X=[0,1]$ and let 
$$
G=([0,1]\times \{1\})\cup(\{0\}\times [0,1]),
$$
 see Figure \ref{figure2}.
\begin{figure}[h!]
	\centering
		\includegraphics[width=20em]{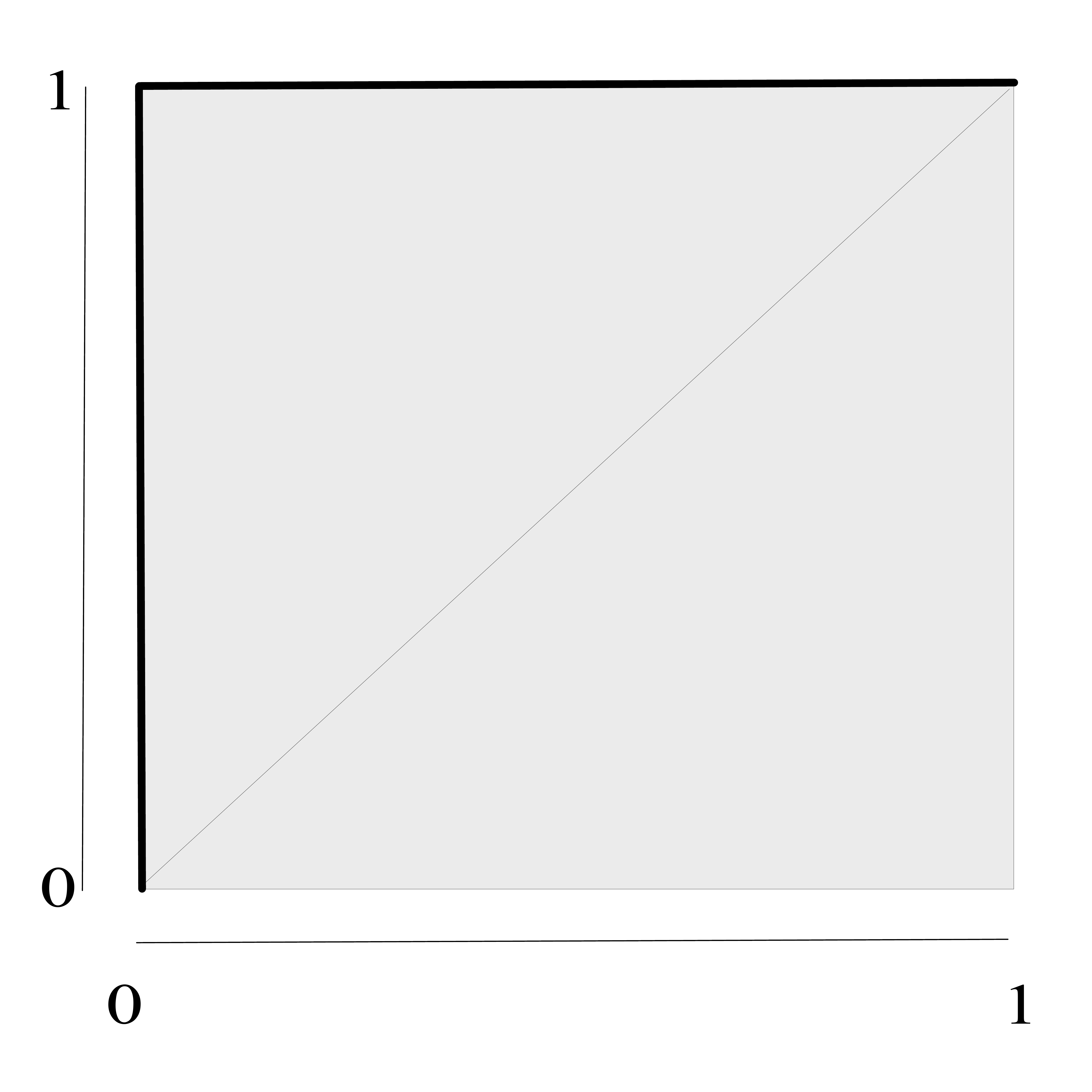}
	\caption{The relation  $G$ from Example \ref{ex2}}
	\label{figure2}
\end{figure} 

\noindent Then $(X,G)$ is a CR-dynamical system. Let $x=0$. Note that for each $\mathbf x\in \tplus_G(x)$, the set $\orbit_G(\mathbf x)$ is not dense in $X$. Also note, that for each $y\in [0,1]$, there is $\mathbf y\in \tplus_G(x)$ such that $\pi_2(\mathbf y)=y$. Therefore, $\orbitt_G(x)$ is dense in $X$. Therefore, 
$$
x\in \trans_3(G)\setminus \trans_2(G).
$$
\end{example}
Note that the relation $G$ in Example \ref{ex2} contains a vertical line that is essential in the example for the  existence of type 3 transitive points in $([0,1],G)$ that are not type 2 transitive.  The following example shows that there are closed relations $G$ on $[0,1]$ such that $G$ contains no vertical lines and there are type 3 transitive points in $([0,1],G)$ that are not type 2 transitive. 
\begin{example}\label{ex4}
Let $X=[0,1]$, let $C$ be the standard ternary Cantor set in $X$,  let $f:X\rightarrow X$ be the Cantor function or the devil's staircase function; see \cite{cantor} for the exact definition and the basic properties of the function  $f$.  Note that $f(C)=[0,1]$.
Let 
$$
G=\left(\left\{\frac{1}{2}\right\}\times C\right)\cup \Gamma(f),
$$
 see Figure \ref{figure4}, where a sketch of an approximation of $G$ is presented.
\begin{figure}[h!]
	\centering
		\includegraphics[width=20em]{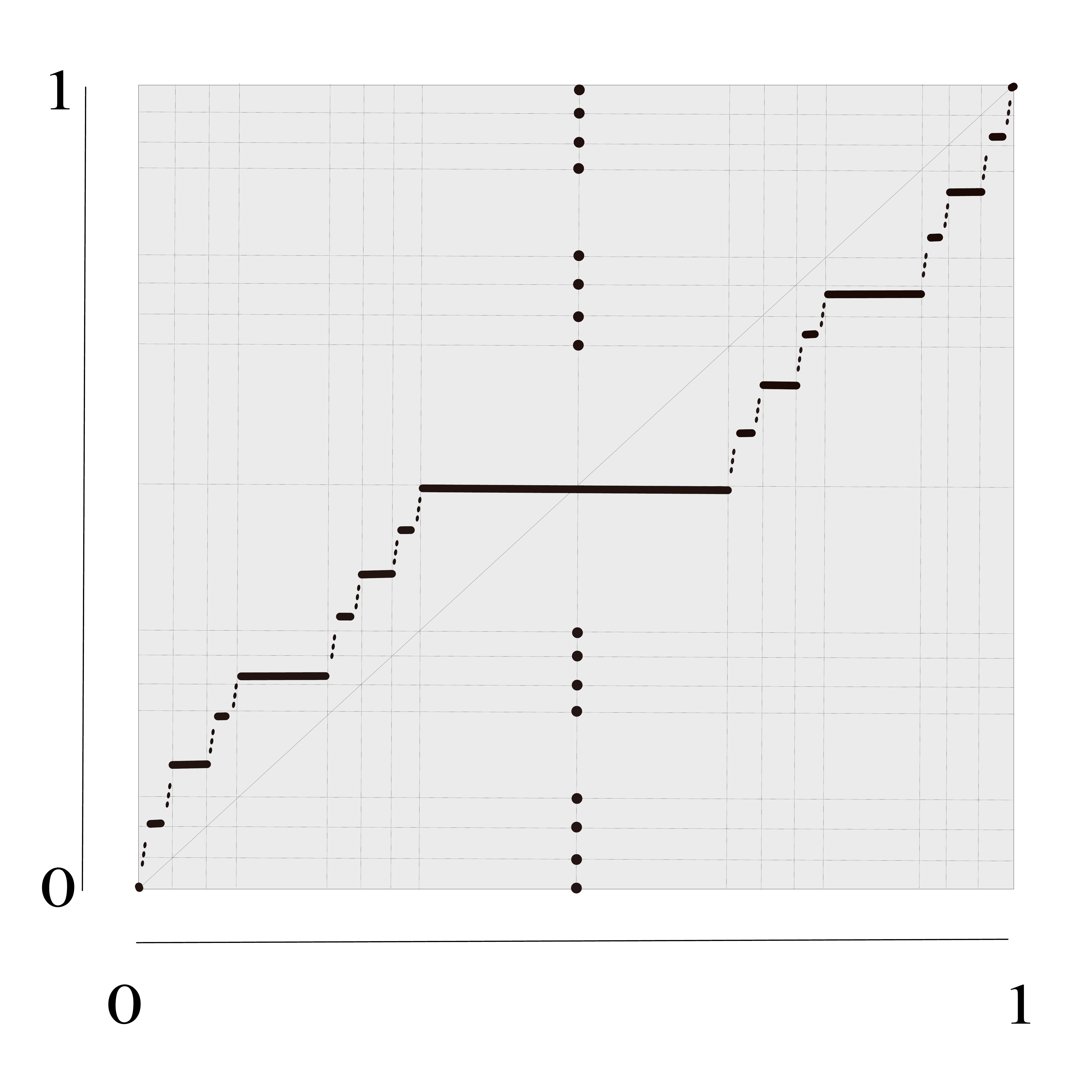}
	\caption{The relation  $G$ from Example \ref{ex4}}
	\label{figure4}
\end{figure} 
\noindent Then $(X,G)$ is a CR-dynamical system. Let $x=\frac{1}{2}$. For each $t\in [0,1]$ let $c_t\in C$ be such a point that $f(c_t)=t$ and let
$$
\mathbf x_t=\Big(\frac{1}{2},c_t,t,f(t),f^2(t),f^3(t),\ldots\Big)
$$
Then 
$$
[0,1]=\bigcup_{t\in [0,1]}\orbit_G(\mathbf x_{t})\subseteq \orbitt_G(x).
$$
It follows that $x\in \trans_3(G)$.  Note that for any $\mathbf x\in \tplus_G(x)$,
$\orbit_G(\mathbf x)$ is not dense in $[0,1]$.  Therefore,  $x\in \trans_3(G)\setminus \trans_2(G)$. 
\end{example}
\begin{definition}
Let $X$ be a compact metric space. We use \emph{\color{blue}$ \isolated(X)$ } to denote the set of isolated points of $X$. 
\end{definition}
\begin{lemma}\label{dense}
Let $X$ be a compact metric space, let $A\subseteq X$ be such that $A$ is not dense in $X$,  and let $x\in X$. If $A\cup \{x\}$ is dense in $X$, then $x\in \isolated(X)$.
\end{lemma}
\begin{proof}
Suppose that $A\cup \{x\}$ is dense in $X$. Then $\Cl(A)\neq X$ while $\Cl(A\cup\{x\})=X$.  Since 
$$
X=\Cl(A\cup\{x\})=\Cl(A)\cup \Cl(\{x\})=\Cl(A)\cup \{x\},
$$
it follows that $\{x\}=X\setminus \Cl(A)$. Therefore,  $x$ is an isolated point in $X$. 
\end{proof}
In the following two examples,  $\isolated(X)\neq \emptyset$ and $\trans_2(G)\neq \emptyset$.  In both examples $\isolated(X)\cap \trans_2(G)\neq \emptyset$.  In Example \ref{fse1},  $\isolated(X)\subseteq  \trans_2(G)$ while in Example \ref{fse2}, $\isolated(X)\not \subseteq  \trans_2(G)$. 
\begin{example}\label{fse1}
Let $X=\{1,2,3\}$ and let $G=\{(1,2),(2,3),(3,1)\}$.  Then $\isolated(X)=\trans_2(G)=\{1,2,3\}$.  Therefore,  $\isolated(X)\neq \emptyset$, $\trans_2(G)\neq \emptyset$, and $\isolated(X)\subseteq  \trans_2(G)$. 
\end{example}
\begin{example}\label{fse2}
Let $X=[0,1]\cup\{2,3\}$ and let  $f:[0,1]\rightarrow [0,1]$ be the tent function defined by $f(t)=2t$ for $t\in [0,\frac{1}{2}]$ and $f(t)=2-2t$ for $t\in [\frac{1}{2},1]$.  Also, let $x\in \tr(f)$, and let 
$$
G=\Gamma(f)\cup \{(2,3),(3,x)\},
$$ 
see Figure \ref{tent}. 
\begin{figure}[h!]
	\centering
		\includegraphics[width=20em]{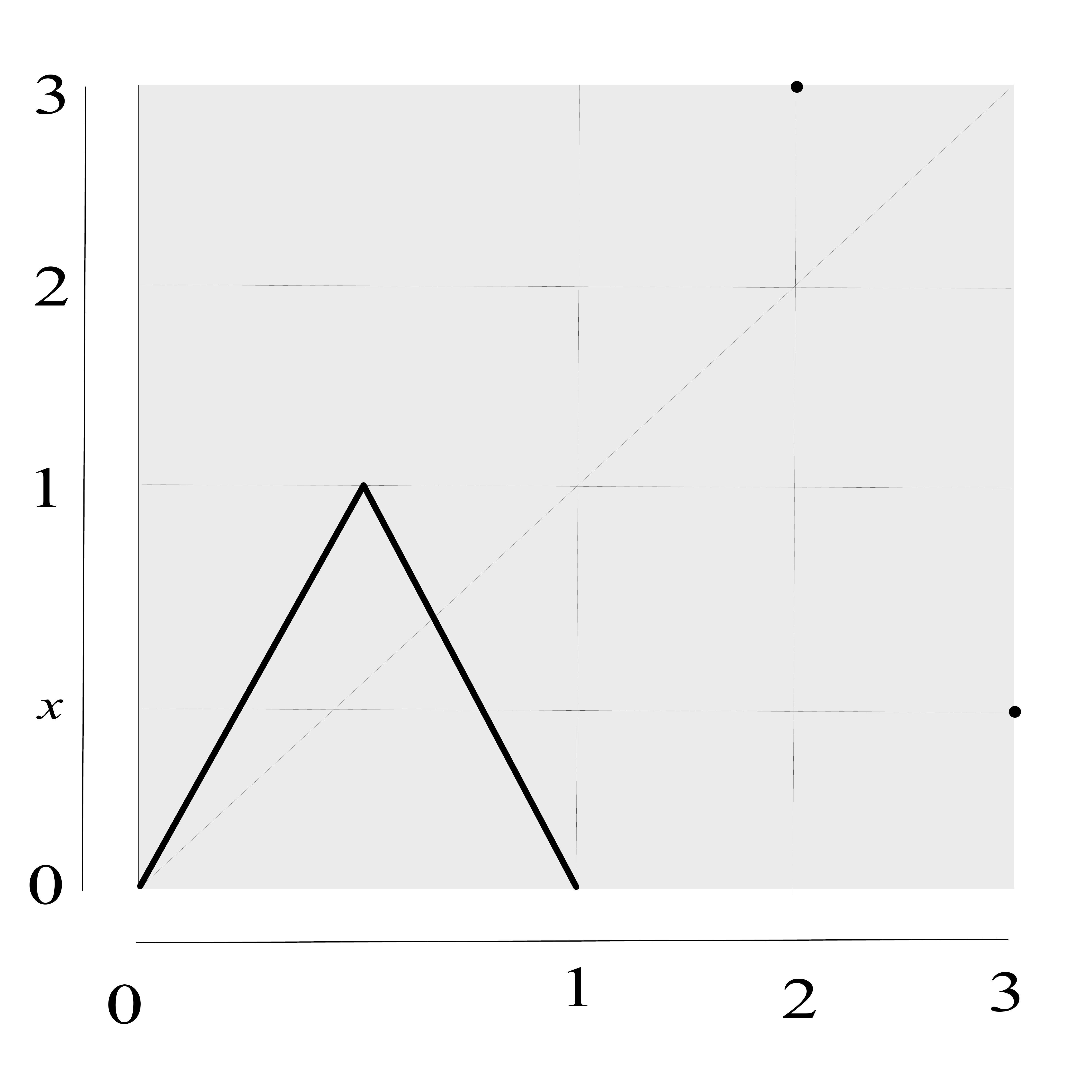}
	\caption{The relation  $G$ from Example \ref{fse2}}
	\label{tent}
\end{figure} 

\noindent Then $\trans_2(G)=\{2\}$ and $\isolated(X)=\{2,3\}$.  It follows that $\isolated(X)\neq \emptyset$,  $\trans_2(G)\neq \emptyset$ and 
$$
 \isolated(X)\not \subseteq  \trans_2(G).
$$
Note that  $\trans_2(G)$ contain exactly one isolated point.
\end{example}
In Example \ref{fse3}, we construct an example of a CR-dynamical system $(X,G)$ such that $X$ is infinite, $\isolated(X)\cap \trans_2(G)\neq \emptyset$,   and $\trans_2(G)$ contains more than just one isolated point.
\begin{example}\label{fse3}
Let $X=[0,1]\cup\{2,3\}$ and let  $f:[0,1]\rightarrow [0,1]$ be the tent function defined by $f(t)=2t$ for $t\in [0,\frac{1}{2}]$ and $f(t)=2-2t$ for $t\in [\frac{1}{2},1]$.  Also, let $x\in \tr(f)$, and let 
$$
G=\Gamma(f)\cup \{(2,3),(3,2),(3,x)\},
$$ 
see Figure \ref{tent1}. 
\begin{figure}[h!]
	\centering
		\includegraphics[width=20em]{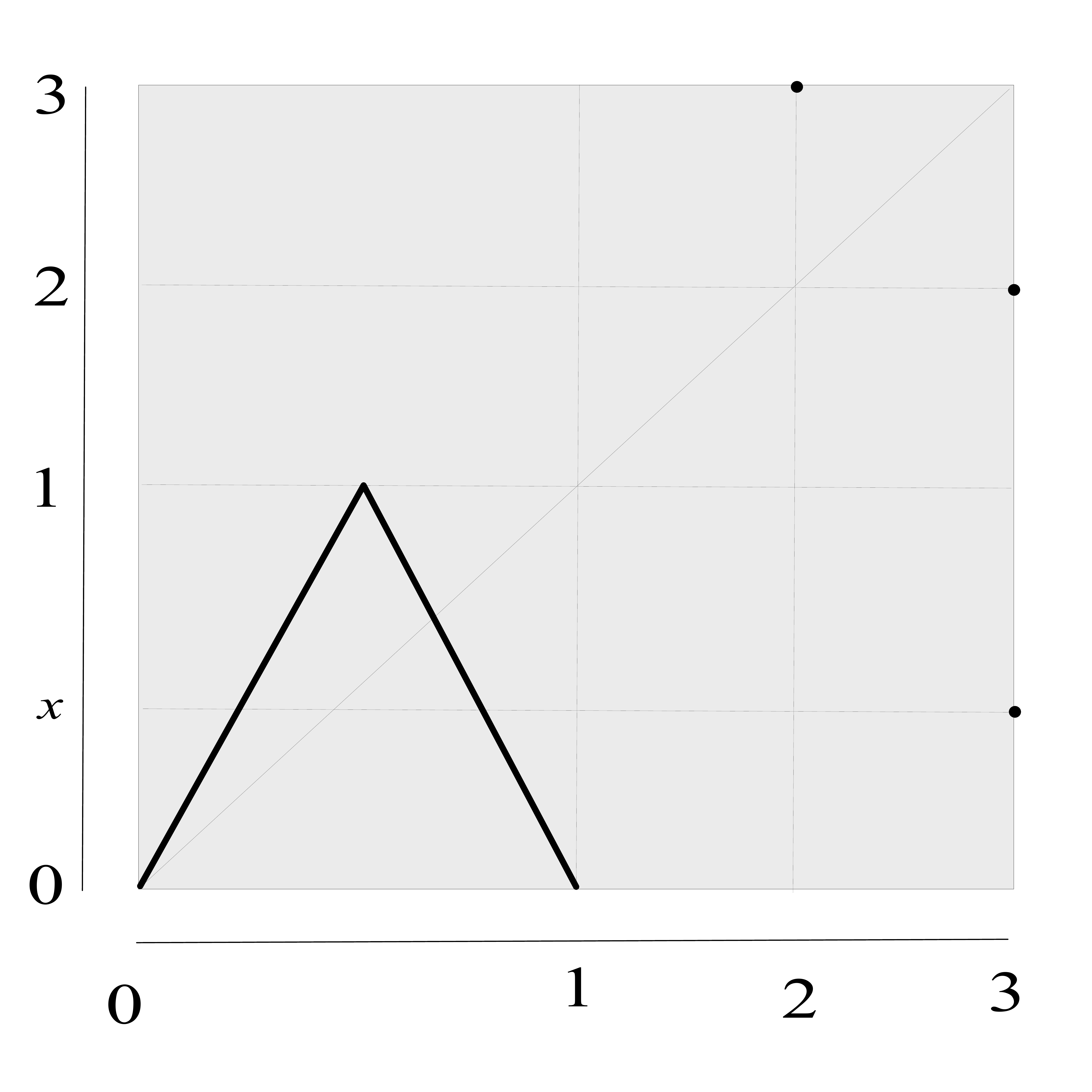}
	\caption{The relation  $G$ from Example \ref{fse3}}
	\label{tent1}
\end{figure} 

\noindent Then $\trans_1(G)=\emptyset$, $\trans_2(G)=\{2,3\}$ and $\isolated(G)=\{2,3\}$.  It follows that $\isolated(X)\cap \trans_2(G)\neq \emptyset$,   and $\trans_2(G)$ contains more than just one isolated point.
\end{example}
In Theorem \ref{Fse11}, we show that if $\isolated(X)\neq \emptyset$, $\trans_1(G)\neq \emptyset$ and $\trans_2(G)\neq \emptyset$, then $\isolated(X)\cap \trans_1(G)\neq \emptyset$ and  $\isolated(X)\cap \trans_2(G)\neq \emptyset$ is always the case.  Also, if $X$ is infinite, then $\trans_1 (G)$ contains exactly one isolated point. 
	\begin{theorem} \label{Fse11}
		Let $(X,G)$ be a CR-dynamical system such that $\isolated(X)\neq \emptyset$. Then the following hold. 
		\begin{enumerate}
			\item If $\trans_1(G) \neq \emptyset$ then $\trans_1 (G)\cap \isolated(X)\neq \emptyset$. In addition, if $X$ is not finite, then  $\trans_1 (G)$ contains exactly one isolated point. 
			\item If $\trans_2(G) \neq \emptyset$ then $\trans_2 (G)\cap \isolated(X)\neq \emptyset$ contains at least one isolated point.
		\end{enumerate}
	\end{theorem}
	\begin{proof}
		Let $k \in \{1,2\}$ and  suppose that  $\trans_k (G) \neq \emptyset$. Let $y \in \trans_k(G)$ and let $\mathbf y \in T_G^+(y)$ such that $\mathcal O_G^{\oplus}(\mathbf y) $ is dense in $X$. Since $\mathcal O_G^{\oplus}(\mathbf y) $ is dense in $X$ all isolated points of $X$ are elements of $\mathcal O_G^{\oplus}(\mathbf y) $. Choose the  smallest $m \in \mathbb N$ such that $\pi_m (\mathbf y )\in \isolated(X)$.  We show that $\pi_m(\mathbf y)\in \trans_1(G)$ by showing that the set $\{ \pi_m (\mathbf y ), \pi_{m+1} (\mathbf y ), \pi_{m+2} (\mathbf y ), \ldots\}$ is dense in $X$.   If $m=1$, then  $\{ \pi_m (\mathbf y ), \pi_{m+1} (\mathbf y ), \pi_{m+2} (\mathbf y ), \ldots\}$ is dense in $X$.  Suppose that $m>1$.  To show that $\{ \pi_m (\mathbf y ), \pi_{m+1} (\mathbf y ), \pi_{m+2} (\mathbf y ), \ldots\}$ is dense in $X$,  let $U$ be any non-empty open set in $X$.  We show that 
		$$
		\{ \pi_m (\mathbf y ), \pi_{m+1} (\mathbf y ), \pi_{m+2} (\mathbf y ), \ldots\}\cap U\neq \emptyset.
		$$
	Since for each $k\in \{1,2,3,\ldots ,m-1\}$, $\pi_k(\mathbf y)\not \in \isolated(X)$, it follows that 
	$$
	U\setminus\{\pi_1(\mathbf y),\pi_2(\mathbf y),\pi_3(\mathbf y),\ldots ,\pi_{m-1}(\mathbf y)\}
	$$	
	is also a non-empty open set in $X$.  Since $\mathcal O_G^{\oplus}(\mathbf y) $ is dense in $X$, it follows that 
	$$
	\mathcal O_G^{\oplus}(\mathbf y)\cap (U\setminus\{\pi_1(\mathbf y),\pi_2(\mathbf y),\pi_3(\mathbf y),\ldots ,\pi_{m-1}(\mathbf y)\})\neq \emptyset.
	$$ 
It follows from $\mathcal O_G^{\oplus}(\mathbf y)=\{\pi_1(\mathbf y),\pi_2(\mathbf y),\pi_3(\mathbf y),\ldots\}$  that  there is a positive integer $k$ such that 
$$
\pi_k(\mathbf y)\in U\setminus\{\pi_1(\mathbf y),\pi_2(\mathbf y),\pi_3(\mathbf y),\ldots ,\pi_{m-1}(\mathbf y)\}.
$$
Choose such a positive integer $k$. Obviously, $k\geq m$.	It follows that 
$$
		\{ \pi_m (\mathbf y ), \pi_{m+1} (\mathbf y ), \pi_{m+2} (\mathbf y ), \ldots\}\cap U\neq \emptyset.
		$$
Therefore,  	$\{ \pi_m (\mathbf y ), \pi_{m+1} (\mathbf y ), \pi_{m+2} (\mathbf y ), \ldots\}$ is dense in $X$.  

For the rest of the proof suppose that $X$ is not finite. We prove that $\trans_1(G)$ contains exactly one of the isolated points. Assume that $x_1, x_2 \in \trans_1(G)\cap \isolated(X)$ such that $x_1\neq x_2$.  Pick any trajectories $\mathbf{x_1} \in T_G^+(x_1)$ and $\mathbf{x_2} \in T_G^+(x_2)$. Since $\mathcal O_G^{\oplus}(\mathbf {x_1}) $ and $\mathcal O_G^{\oplus}(\mathbf {x_2}) $ are dense in $X$ there are positive integers $m$ and $n$ such that $\pi_m (\mathbf {x_1})=x_2$ and $\pi_n (\mathbf {x_2})=x_1$.  Note that $m,n>1$ since $x_1\neq x_2$.  Also, note that $\pi_m (\mathbf {x_1})=\pi_1(\mathbf {x_2})$ and $\pi_n (\mathbf {x_2})=\pi_1(\mathbf {x_1})$. Let $\mathbf z\in \prod_{i=1}^{\infty}X$ be defined by
$$
\mathbf z=(\pi_1(\mathbf {x_1}), \pi_2(\mathbf {x_1}), \ldots,\pi_m(\mathbf {x_1}), \pi_2(\mathbf{x_2}), \pi_3(\mathbf{x_2}), \ldots,\pi_n(\mathbf {x_2}),  \pi_2(\mathbf{x_1}),  \pi_3(\mathbf{x_1}),\ldots), 
$$
Then $\mathbf z\in T_G^+(x_1)$ and 
$$
\Orbit(\mathbf z)=\{\pi_1(\mathbf {x_1}), \pi_2(\mathbf {x_1}), \ldots,\pi_m(\mathbf {x_1}), \pi_2(\mathbf{x_2}), \pi_3(\mathbf{x_2}), \ldots,\pi_{n-1}(\mathbf {x_2})\}.
$$
It follows that $\Orbit(\mathbf z)$ is finite and since $X$ is infinite,  $\Orbit(\mathbf z)$ is not dense in $X$,  which is a contradiction since $x_1\in \trans_1(G)$ and $\mathbf z\in T_G^+(x_1)$.
	\end{proof}
	In the following example we show that it may happen that  $\trans_3(G)\neq \emptyset$ and  $\trans_3(G)$ contains no isolated points (even in the case when $\isolated (X)\neq \emptyset$).
	\begin{example}\label{ff}
		Let $X=[0,1] \cup \{2\}$, let  $y\in (0,1]$ and let 
		$$
		G=([0,1]\times \{1\})\cup(\{0\}\times [0,1]) \cup \{(0,2), (2,y) \},
		$$
		see Figure \ref{tent11}. 
\begin{figure}[h!]
	\centering
		\includegraphics[width=20em]{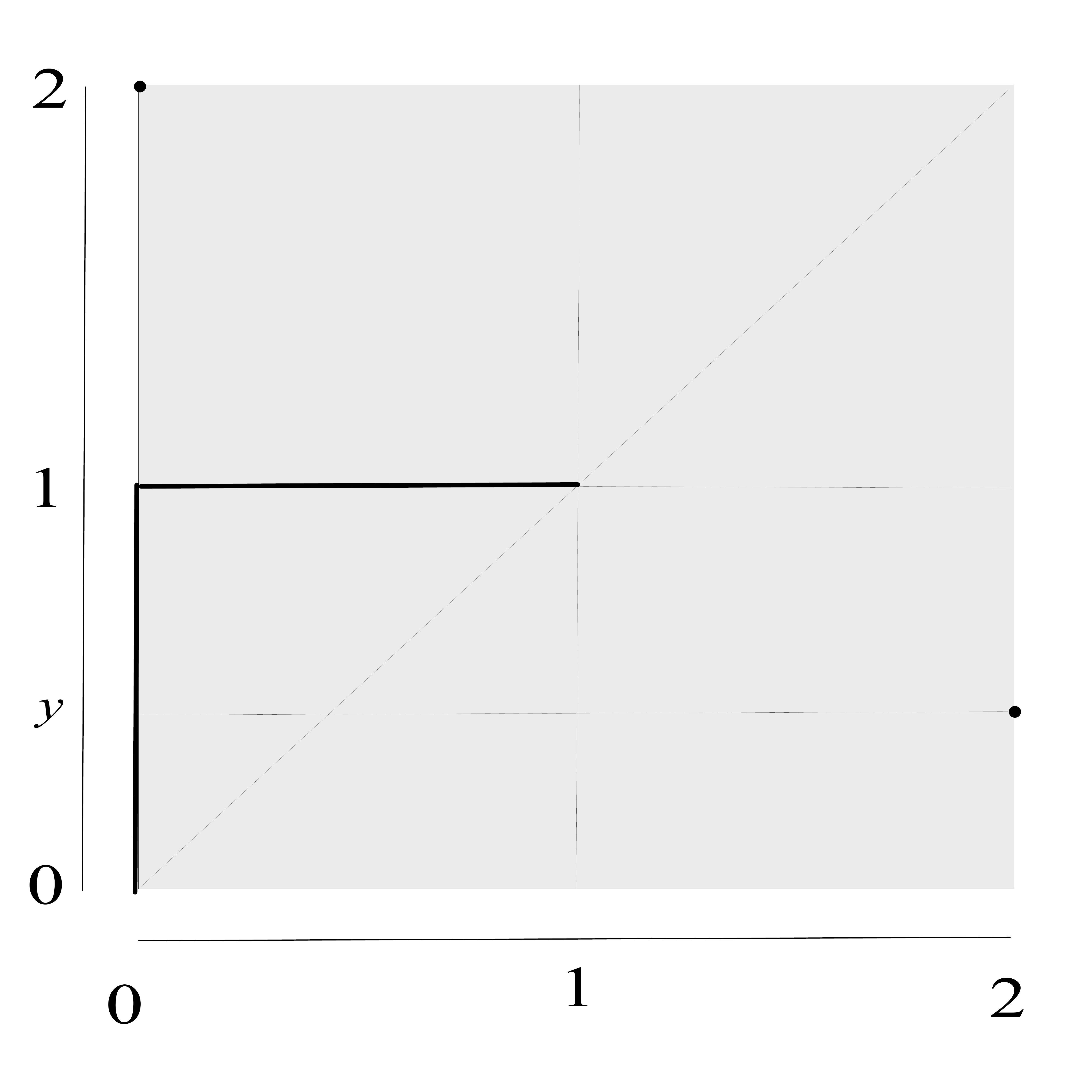}
	\caption{The relation  $G$ from Example \ref{ff}}
	\label{tent11}
\end{figure} Note that  $\isolated(X)=\{2\}$.  Also,  note that $$\orbitt_G(2)=\orbit_G(2)=\{2,y,1\}$$ 
		and, therefore,  $\orbitt_G(2)$ is not dense in $X$. So, $2 \notin \trans_3(G)$.   Next, let $x\in (0,1]$. Then $$\orbitt_G(x)=\orbit_G(x)=\{x,1\}$$ which is not a dense subset of $X$.  Therefore,  $x \notin \trans_3(G)$. It follows that $ \trans_3(G)=\{0\}$.  Therefore,
		$\trans_3(G)\cap \isolated(X)=\emptyset$.
	\end{example}

Note that there are examples of CR-dynamical systems $(X,G)$ such that $X$ is infinite,  with an isolated point $x$ such that for some $y\in X$,  $(x,y)\in G$ and $y\in \trans_2(G)$.  The following example demonstrates this.  In Theorem \ref{sinaG1}, we show that for $1$-transitive points, this is not the case.

\begin{example}\label{tistile0}
		Let $X=[0,1]\cup\{2\}$ and let $f:[0,1]\rightarrow [0,1]$ be the tent-map defined by $f(t)=2t$ for $t\leq \frac{1}{2}$ and $f(t)=2-2t$ for $t\geq \frac{1}{2}$,  let $t_0\in \tr(f)$ and let  $$G=\Gamma(f)\cup \{(2,t_0),(t_0,2)\},$$  see Figure \ref{figgure0}.
		\begin{figure}[h!]
			\centering
			\includegraphics[width=20em]{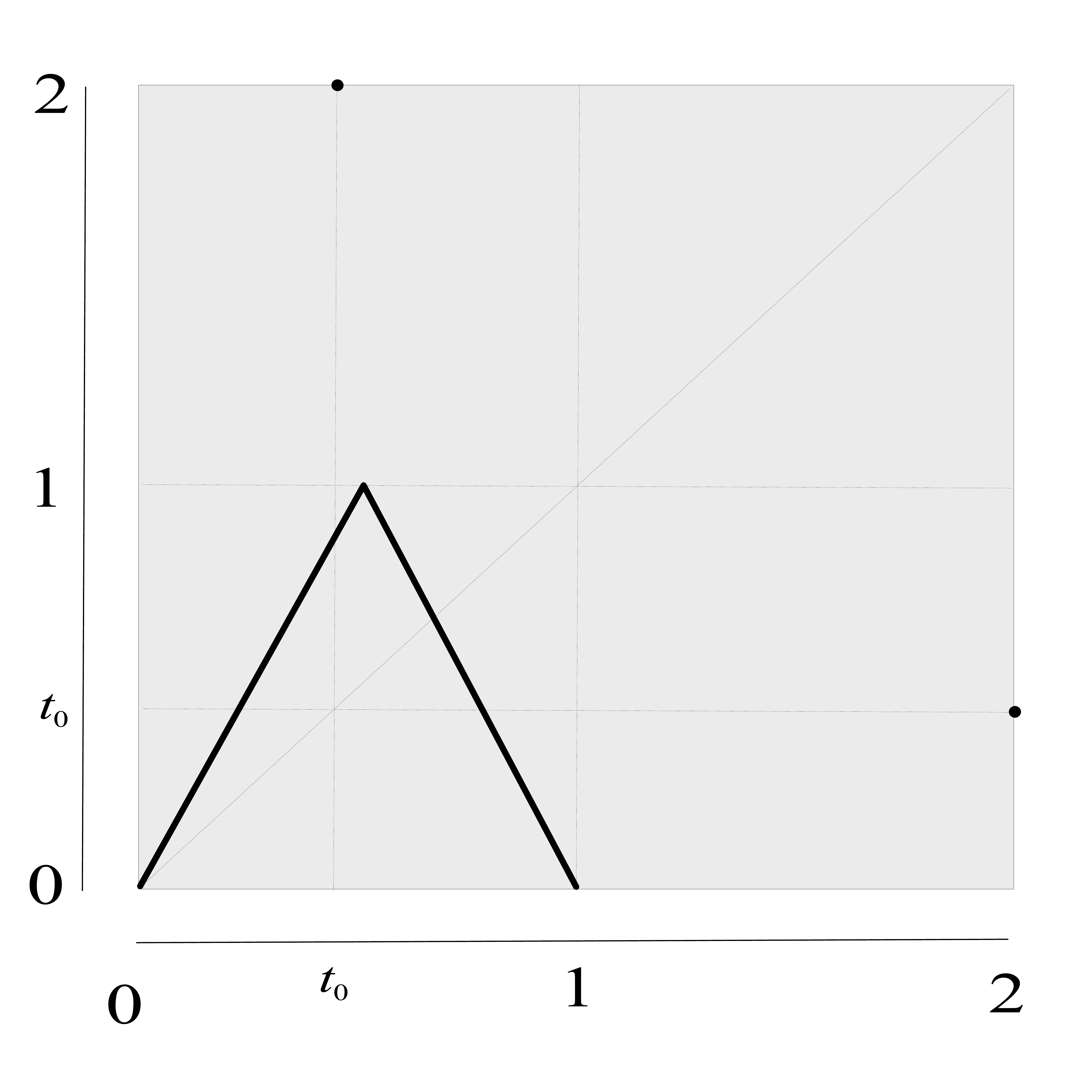}
			\caption{The relation  $G$ from Example \ref{tistile0}} 
			\label{figgure0}
		\end{figure} 
Note that for $x=2$ and $y=t_0$, 	
$$
x\in \isolated(X), (x,y)\in G  \text{ and } y\in \trans_2(G).
$$
	\end{example}
\begin{theorem}\label{sinaG1}
 Let $(X,G)$ be a CR-dynamical system such that $X$ is infinite,   let $x\in \isolated(X)$  and let $y\in X$ be any point.  
\begin{enumerate}
\item If $y\in \trans_2(G)$, then either $x=y$ or there is a positive integer $n$ and a point $\mathbf x\in \star_{i=1}^nG$ such that  $\pi_1(\mathbf x)=y$ and $\pi_{n+1}(\mathbf x)=x$.
\item If $(x,y)\in G$, then $y\not\in \trans_1(G)$.
\end{enumerate}
\end{theorem}
\begin{proof}
First, we prove the first statement.  If $x=y$, there is nothing to prove.  Suppose that $x\neq y$.  Let $y\in \trans_2(G)$ and let $\mathbf y\in \star_{i=1}^{\infty}G$ be a trajectory of $y$ such that the forward orbit  $\orbit_G(\mathbf y)$ is dense in $X$.  Since $x\in \isolated(X)$, $\{x\}$ is open in $X$ and, therefore, $\orbit_G(\mathbf y)\cap \{x\}\neq \emptyset$.  It follows that there is a positive integer  $m$  such that $\pi_m(\mathbf y)=x$.  Let $m$ be such a positive integer and let $n=m-1$ (note that $n>0$).  Then let 
$$
\mathbf x=(\pi_1(\mathbf y),\pi_2(\mathbf y),\pi_3(\mathbf y),\ldots,\pi_{n+1}(\mathbf y)).
$$
Obviously, $\pi_1(\mathbf x)=y$ and $\pi_{n+1}(\mathbf x)=x$.

To prove the second claim, suppose that $y\in \trans_1(G)$. It follows that $y\in \trans_2(G)$ and we have just proved that either $x=y$ or there is a positive integer $n$ and a point $\mathbf x\in \star_{i=1}^nG$ such that  $\pi_1(\mathbf x)=y$ and $\pi_{n+1}(\mathbf x)=x$.  Next, we observe each of these two possibilities.
\begin{enumerate}
\item If $x=y$, then $\mathbf y=(y,y,y,\ldots)\in \tplus_G(y)$
 such that $\orbit_G(\mathbf y)=\{y\}$ is not dense in $X$ since $X$ is not degenerate -- a contradiction.
 \item Let $n$ be a positive integer and let $\mathbf x\in \star_{i=1}^nG$ be such a point that  $\pi_1(\mathbf x)=y$ and $\pi_{n+1}(\mathbf x)=x$.  Let
 $$
 \mathbf y=(\pi_1(\mathbf x),\pi_2(\mathbf x),\pi_3(\mathbf x),\ldots,\pi_{n+1}(\mathbf x),\pi_1(\mathbf x),\pi_2(\mathbf x),\pi_3(\mathbf x),\ldots,\pi_{n+1}(\mathbf x),\pi_1(\mathbf x),\ldots).
 $$
 Then $\orbit_G(\mathbf y)=\{\pi_1(\mathbf x),\pi_2(\mathbf x),\pi_3(\mathbf x),\ldots,\pi_{n+1}(\mathbf x)\}$ and since $X$ is infinite,  the orbit $\orbit_G(\mathbf y)$ of $y$ is not dense in $X$ -- a contradiction. 
\end{enumerate}
Therefore,  $y\not\in \trans_1(G)$. 
\end{proof}

\begin{theorem}\label{319}
Let $(X,G)$ be a CR-dynamical system and let $x\in \legal(G)$. Then the following hold.
\begin{enumerate}
\item\label{1} If $x\in \intrans(G)$, then for each $y\in X$,
$$
(x,y)\in G \Longrightarrow  y\in \intrans(G).
$$
\item\label{2} If $x\in \intrans(G)$, then for each $\mathbf x\in \tplus_G(x)$ and for each positive integer $n$, 
$$
\pi_n(\mathbf x)\in \intrans(G).
$$
\item\label{3} If $x$ is  not an isolated point and if $x\in \trans_1(G)$, then for each $y\in \legal(G)$,
$$
(x,y)\in G \Longrightarrow  y\in \trans_1(G).
$$
\item\label{4} If $X$ has  no isolated points and if $x\in \trans_1(G)$, then for each $\mathbf x\in \tplus_G(x)$ and for each positive integer $n$,
$$
\pi_n(\mathbf x)\in \trans_1(G).
$$
\item\label{5} If $x$ is  not an isolated point and if $x\in \trans_2(G)$, then there is $y\in \legal(G)$ such that
$$
(x,y)\in G \textup{ and }  y\in \trans_2(G).
$$
\item\label{6} If $X$ has  no isolated points  and if $x\in \trans_2(G)$, then there is $\mathbf x\in \tplus_G(x)$ such that $\orbit_G(\mathbf x)$ is dense in $X$ and for each positive integer $n$,
$$
\pi_n(\mathbf x)\in \trans_2(G).
$$
\end{enumerate}
\end{theorem}

\begin{proof}
To prove \ref{1}, let $x\in \intrans(G)$ and let $y\in X$ such that $(x,y)\in G$. If $y\not \in \intrans(G)$, then $y\in \trans_3(G)$. It follows that $\orbitt_G(y)$ is dense in $X$. Since $\orbitt_G(y)\subseteq \orbitt_G(x)$, it follows that $\orbitt_G(x)$ is dense in $X$, meaning that $x\in \trans_3(G)$  -- a  contradiction.

To prove \ref{2}, let $x\in \intrans(G)$ and let $\mathbf x\in \tplus_G(x)$. Suppose that $n$ is a positive integer such that $\pi_n(\mathbf x)\not \in \intrans(G)$. It follows that $\orbitt_G(\pi_n(\mathbf x))$ is dense in $X$.   Since $\orbitt_G(\pi_n(\mathbf x))\subseteq \orbitt_G(x)$, it follows that $\orbitt_G(x)$ is dense in $X$.  Therefore,  $x\not \in \intrans(G)$  -- a  contradiction.  It follows that for each positive integer $n$,  $\pi_n(\mathbf x)\in \intrans(G)$.

To prove \ref{3}, suppose that $x$ is not an isolated point in $X$ and that $x\in \trans_1(G)$ and let $y\in \legal(G)$ be such that $(x,y)\in G$.  Suppose that $y\not \in \trans_1(G)$. Then there is $\mathbf y\in \tplus_G(y)$ such that $\orbit_G(\mathbf y)$ is not dense in $X$.  Let 
$$
\mathbf x=(x,y,\pi_2(\mathbf y), \pi_3(\mathbf y), \pi_4(\mathbf y), \ldots).
$$
Then $\mathbf x\in \tplus_G(x)$ is such a point that $\Cl(\orbit_G(\mathbf x))=\Cl(\orbit_G(\mathbf y))\cup\{x\}$. Therefore,  by Lemma \ref{dense}, $\orbit_G(\mathbf x)$ is not dense in $X$  -- a  contradiction. 

To prove \ref{4},  suppose that $X$ has no isolated points and that $x\in \trans_1(G)$ and let $\mathbf x\in \tplus_G(x)$. Suppose that there is a positive integer $n$ such that $\pi_n(\mathbf x)\not \in \trans_1(G)$. Let $n$ be the smallest such integer.  Note that $n>1$.  Let $t=\pi_{n-1}(\mathbf x)$ and $y=\pi_{n}(\mathbf x)$.  Then $t\in \trans_1(G)$, $y\in \legal(G)$,  $(t,y)\in G$, and $y\not \in \trans_1(G)$.  This contradicts \ref{3}.  

To prove \ref{5}, suppose that $x$ is not an isolated point in $X$ and that $x\in \trans_2(G)$. We show that there is $y\in \legal(G)$ such that 
$$
(x,y)\in G \textup{ and }  y\in \trans_2(G).
$$
Since $x\in \trans_2(G)$, there is $\mathbf x\in \tplus_G(x)$ such that $\orbit_G(\mathbf x)$ is dense in $X$.  Let $\mathbf x$ be such a point and let $y=\pi_2(\mathbf x)$.  Then  $(x,y)\in G$. Suppose that $y\not \in \trans_2(G)$ and let 
$$
\mathbf y=(\pi_2(\mathbf x), \pi_3(\mathbf x), \pi_4(\mathbf x), \ldots).
$$
Then $\orbit_G(\mathbf x)=\orbit_G(\mathbf y)\cup\{x\}$. Since $\orbit_G(\mathbf x)$ is dense in $X$ while $\orbit_G(\mathbf y)$ is not, it follows from Lemma \ref{dense} that $x$ is an isolated point in $X$  -- a  contradiction.   

Finally, we prove \ref{6}.  Suppose that $X$ has no isolated points and suppose that $x\in \trans_2(G)$.  Since $x\in \trans_2(G)$, there is $\mathbf x\in \tplus_G(x)$ such that $\orbit_G(\mathbf x)$ is dense in $X$.  Let $\mathbf x$ be such a point and let for each positive integer $n$,
$$
\mathbf y_n=(\pi_n(\mathbf x), \pi_{n+1}(\mathbf x), \pi_{n+2}(\mathbf x), \ldots).
$$
  We show that for each positive integer $n$,  $\pi_n(\mathbf x)\in \trans_2(G)$ by showing that for each positive integer $n$,   $\orbit_G(\mathbf y_n)$ is dense in $X$.  Suppose that there is a positive integer $n$ such that   $\orbit_G(\mathbf y_n)$ is not dense in $X$.  Let $n$ be the smallest such integer.  Note that $n>1$.  Then $\orbit_G(\mathbf y_{n-1})$ is dense in $X$ and $\orbit_G(\mathbf y_{n})$ is not dense in $X$. Since $\orbit_G(\mathbf y_{n-1})=\orbit_G(\mathbf y_{n})\cup\{\pi_{n-1}(\mathbf x)\}$, it follows from Lemma \ref{dense} that $\pi_{n-1}(\mathbf x)$ is an isolated point in $X$  -- a  contradiction.  This completes the proof. 
\end{proof}
In the following example, we show that \ref{3} from Theorem \ref{319} cannot be generalized to the case when $x\in \isolated(X)$. 
\begin{example}\label{dens}
Let $X=\{0,1\}$ and let $G=\{(0,1),(1,1)\}$. Let $x=0$ and $y=1$. Note that  $(x,y)\in G$ and that $x$ is an isolated point of $X$ and that $x\in \trans_1(G)$ while $y\not\in \trans_1(G)$.
\end{example}
In Example \ref{ex3}, we demonstrate why Theorem \ref{319} does not deal with the set  $\trans_3(G)$.
\begin{example}\label{ex3}
Let $X=[0,1]$ and let 
$
G=([0,1]\times \{1\})\cup(\{0\}\times [0,1]),
$
 see Figure \ref{figure3}.
\begin{figure}[h!]
	\centering
		\includegraphics[width=20em]{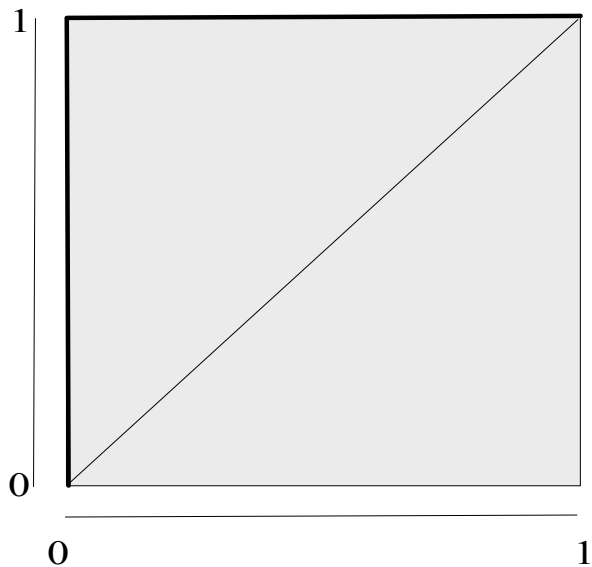}
	\caption{The relation  $G$ from Example \ref{ex3}}
	\label{figure3}
\end{figure} 
We proved in Example \ref{ex2} that $0\in \trans_3(G)$. Note that for each $y\in (0,1]$, $(0,y)\in G$ while $y\not \in \trans_3(G)$.
\end{example}
\begin{theorem}\label{transdense}
Let $(X,G)$ be a CR-dynamical system and let $k\in \{1,2\}$.  If $\isolated(X)=\emptyset$ and $\trans_k(G)\neq \emptyset$, then $\trans_k(G)$  is  dense in $X$.
\end{theorem}
\begin{proof}
Let $k=1$ and suppose that $X$ does not have isolated points and that $\trans_1(G)\neq \emptyset$. Let $x\in \trans_1(G)$ and let $\mathbf x\in \tplus_G(x)$. Then $\Cl(\orbit_G(\mathbf x))=X$. By \ref{4} from Theorem \ref{319}, for each positive integer $n$,
$$
\pi_n(\mathbf x)\in \trans_1(G).
$$
Therefore, $\orbit_G(\mathbf x)\subseteq \trans_1(G)$ and $\Cl(\trans_1(G))=X$ follows. 

Let $k=2$ and suppose that $X$ does not have isolated points and that $\trans_2(G)\neq \emptyset$.  Let $x\in \trans_2(G)$ and let $\mathbf x\in \tplus_G(x)$ such that $\Cl(\orbit_G(\mathbf x))=X$ and such that for each positive integer $n$,
$$
\pi_n(\mathbf x)\in \trans_2(G).
$$
Such a point does exist by Theorem \ref{319}.
Therefore, $\orbit_G(\mathbf x)\subseteq \trans_2(G)$ and $\Cl(\trans_2(G))=X$ follows. 
\end{proof}
\begin{observation}
In Example \ref{dens}, a CR-dynamical system is given such that $\trans_1(G)\neq \emptyset$ and $\trans_2(G)\neq \emptyset$  while neither $\trans_1(G)$ or $\trans_2(G)$ is dense in $X$.  Note that in this example $\isolated{(X)}\neq \emptyset$.
\end{observation}

 The following corollary easily follows from  Theorem \ref{transdense}.
	\begin{corollary}
		Let $(X,G)$ be a CR-dynamical system and let $k\in \{1,2\}$.  If $\isolated(X)=\emptyset$ and $ \trans_k(G)\neq \emptyset$, then $ \illegal(G)=\emptyset$.
	\end{corollary}
	\begin{proof}
		Assume that $\illegal(G) \neq \emptyset.$ By Theorem \ref{illegal_open}, the set $\illegal(G)$ is an open set in $X$. Choose any point $x\in \trans_k(G)$ and let $\mathbf x \in T_G^+ (x)$ be  such that $\orbit_G (\mathbf x)$ is dense in $X$.  Since $\illegal(G)$ is open in $X$ and $\orbit_G (\mathbf x)$ is dense in $X$, it follows that 
		$$
		\illegal(G)\cap \orbit_G (\mathbf x)\neq \emptyset.
		$$  
Since $\orbit_G (\mathbf x)=\{\pi_1(\mathbf x),\pi_2(\mathbf x),\pi_3(\mathbf x),\ldots\}$,  there exists $m \in \mathbb N$ such that 
$$
\pi_m(\mathbf x) \in \illegal(G).
$$
 But notice that $T_G^+ (\pi_m(\mathbf x)) \neq \emptyset$. Hence, $\pi_m(\mathbf x) \in \legal(G)$ and we get a contradiction.
	\end{proof}

Note that there are CR-dynamical systems $(X,G)$ such that $\trans_3(G)$ is non-empty and $\trans_3(G)$ is not dense in $X$; Example \ref{ex3} is such an example.  Next, we give more examples of CR-dynamical systems such that $\trans_3(G)\neq \emptyset$ and $\trans_3(G)$ is not dense in $X$. Example \ref{ex31} does not have isolated points, Example \ref{ex32} does have isolated points.  In addition to that, we show in both examples that for each $x\in \trans_3(G)$,  $x\not \in \trans_2(G)$.  This serves as a motivation for Definition \ref{666}.
\begin{example}\label{ex31}
Let $X=[0,1]$, and let 
$$
f_{1}:\Big[0,\frac{1}{2}\Big]\rightarrow \Big[0,\frac{1}{2}\Big]
$$
 be defined by $f_1(t)=2t$ if $t\in[0,\frac{1}{4}]$ and $f_1(t)=1-2t$ if $t\in[\frac{1}{4},\frac{1}{2}]$ and let 
 $$
 f_{2}:\Big[\frac{1}{2},1\Big]\rightarrow \Big[\frac{1}{2},1\Big]
 $$
  be defined by $f_2(t)=2t-\frac{1}{2}$ if $t\in[\frac{1}{2},\frac{3}{4}]$ and $f_2(t)=\frac{5}{2}-2t$ if $t\in[\frac{3}{4},1]$.   Let $x_1\in \tr(f_1)$ and $x_2\in \tr(f_2)$, and let  $G$ be defined as follows:
$$
G=\Gamma(f_1)\cup \Gamma(f_2)\cup\{(0,x_1),(0,x_2)\},
$$
see Figure \ref{fig31}.
\begin{figure}[h!]
	\centering
		\includegraphics[width=20em]{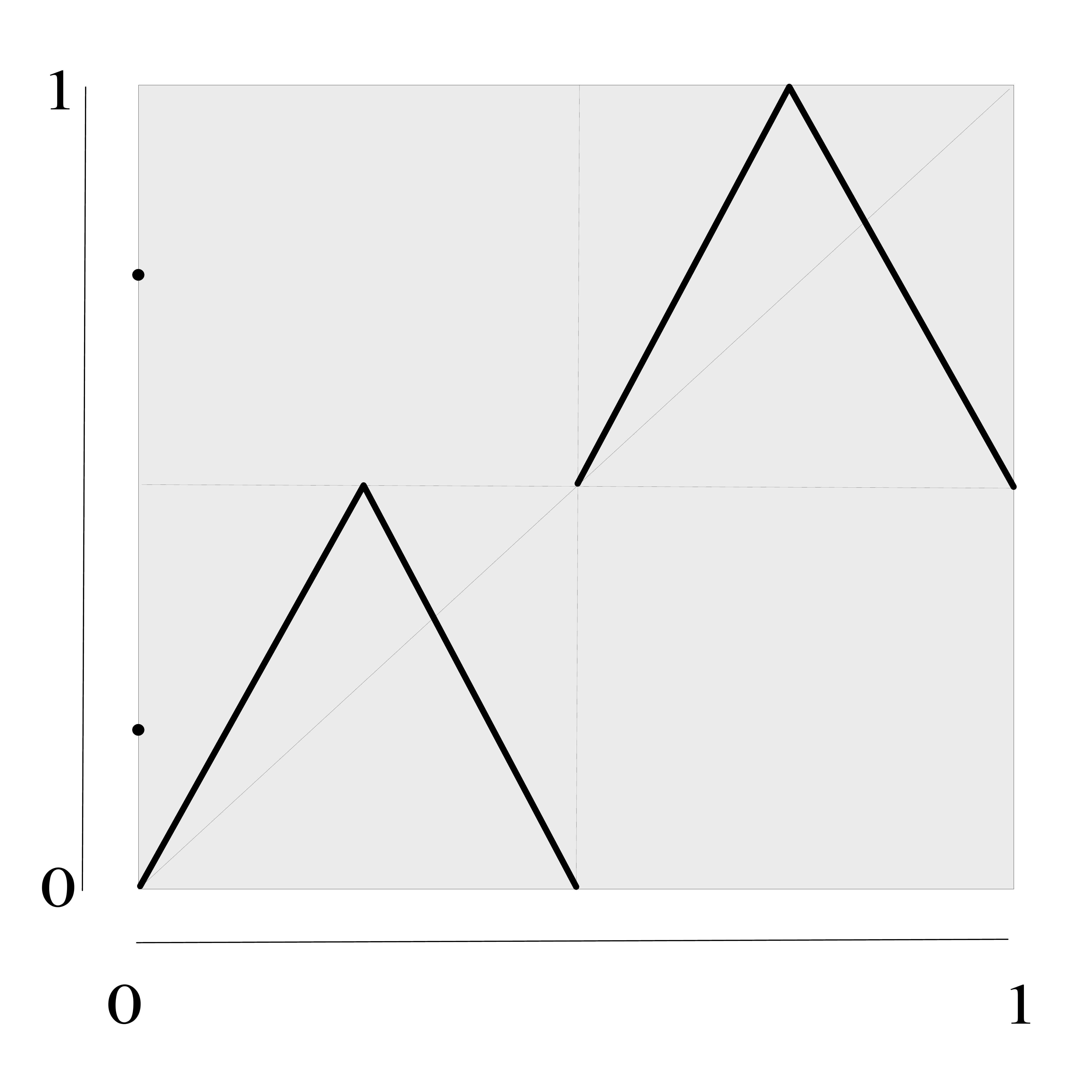}
	\caption{The relation  $G$ from Example \ref{ex31}}
	\label{fig31}
\end{figure} 

\noindent Let $x=0$. We show first that $x$ is a type 3 transitive point in $([0,1],G)$ and is not a type 2 transitive point in $([0,1],G)$. To see that $x\not \in \trans_2(G)$, let $\mathbf x\in \tplus_G(x)$ be any point.  Then there are the following possible cases.
\begin{enumerate}
\item $\mathbf x=(0,0,0,\ldots)$. In this case $\orbit_G(\mathbf x)=\{0\}$ and, therefore, $\orbit_G(\mathbf x)$ is not dense in $X$. 
\item There is a positive integer $n$ such that $\pi_n(\mathbf x)=x_1$.  In this case,  $\orbit_G(\mathbf x)\subseteq [0,\frac{1}{2}]$ and, therefore,  $\orbit_G(\mathbf x)$ is not dense in $X$. 
\item There is a positive integer $n$ such that $\pi_n(\mathbf x)=x_2$. In this case,  $\orbit_G(\mathbf x)\subseteq [\frac{1}{2},1]$ and, therefore,  $\orbit_G(\mathbf x)$ is not dense in $X$. 
\end{enumerate}
Next, let  $\mathbf x_1,\mathbf x_2\in \tplus_G(x)$ be such that $\pi_2(\mathbf x_1)=x_1$ and $\pi_2(\mathbf x_2)=x_2$.  Then $\orbit_G(\mathbf x_1)\cup \orbit_G(\mathbf x_2)$ is dense in $X$. Therefore, $x\in \trans_3(G)$. 
Note that $\trans_3(G)=\{0\}$ and, therefore, $\trans_3(G)$ is not dense in $X$. 
\end{example}
\begin{example}\label{ex32}
For each positive integer $n$, let 
$$
a_1^n,a_2^n,a_3^n,\ldots ,a_{n}^n\in \Big(\frac{2^{n+1}-1}{2^{n+1}},\frac{2^{n+2}-1}{2^{n+2}}\Big)
$$ 
be $n$  distinct points and let 
$$
A_n=\{a_1^n,a_2^n,a_3^n,\ldots ,a_{n}^n\} 
$$
and
$$
 G_n=\left\{\frac{2^{n+1}-1}{2^{n+1}}\right\}\times A_n
$$ 
Let 
$$
X=\{1\}\cup \left\{\frac{2^{n-1}-1}{2^{n-1}} \ | \ n \textup{ is a positive integer}\right\}\cup \bigcup_{n=1}^{\infty}A_n
$$
 and let 
 $$
G=\{(x,x)  \  |  \  x\in X \}\cup \left\{\left(0,\frac{2^{n-1}-1}{2^{n-1}}\right) \ | \ n \textup{ is a positive integer}\right\}\cup \bigcup_{n=1}^{\infty}G_n,
$$
see Figure \ref{fig32}.
\begin{figure}[h!]
	\centering
		\includegraphics[width=20em]{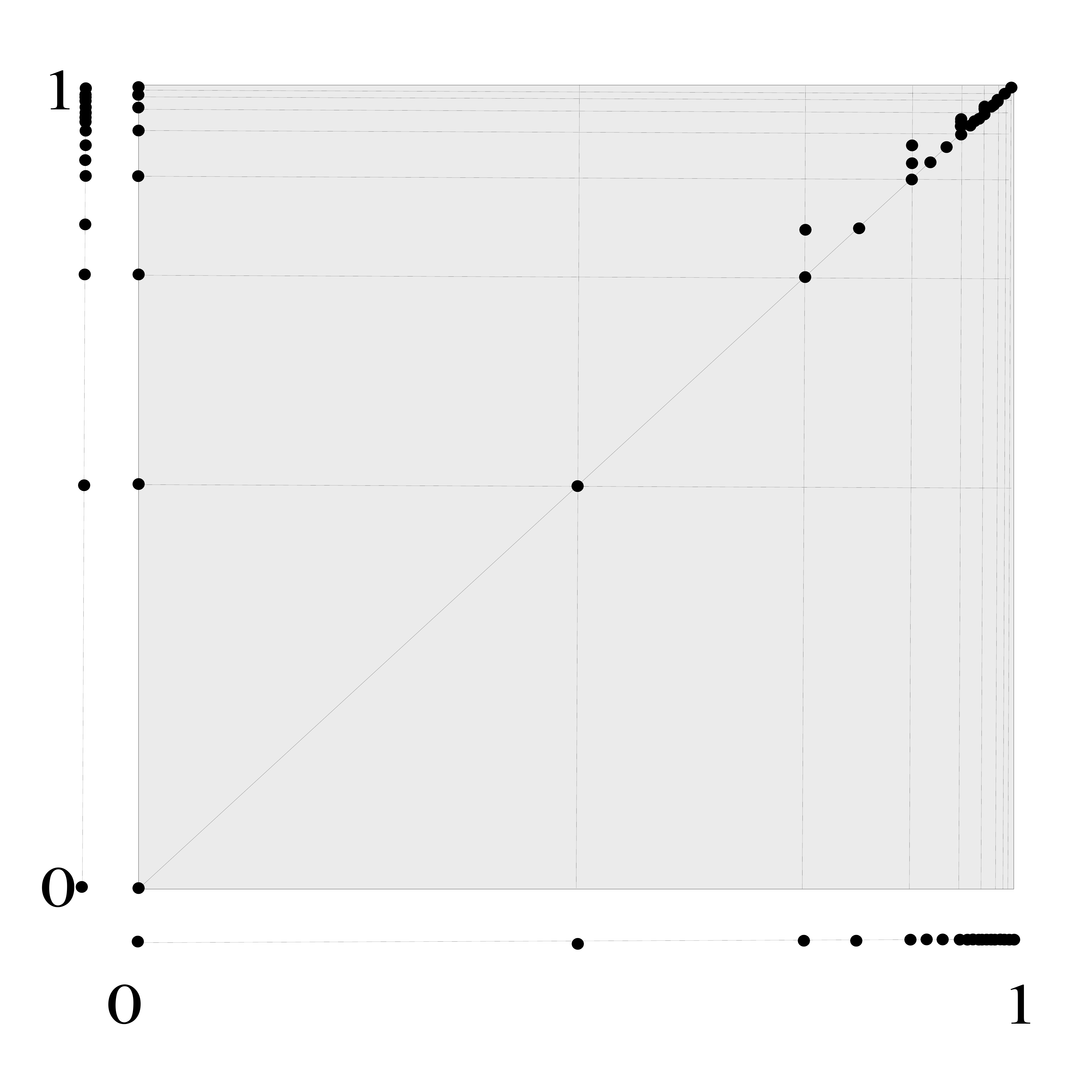}
	\caption{The relation  $G$ from Example \ref{ex32}}
	\label{fig32}
\end{figure} 

\noindent Then $0\in \trans_3(G)$ but $0\not \in \trans_2(G)$.  To see this, let $\mathbf x\in \tplus_G(0)$.  If $\pi_k(\mathbf x)=0$ for each positive integer $k$, then $\mathbf x=(0,0,0,\ldots)$ and $\orbit_G(\mathbf x)=\{0\}$, which is not dense in $X$.  Suppose that there is a positive integer $k$ such that $\pi_k(\mathbf x)\neq 0$ and let 
$$
m=\min\{k \ | \ k \textup{ is a positive integer such that } \pi_k(\mathbf x)\neq 0\}.
$$
Note that $m\geq 2$.  Then there is a non-negative integer $n$ such that 
$$
\pi_m(\mathbf x)=\frac{2^{n+1}-1}{2^{n+1}} \textup{ or } \pi_m(\mathbf x)=1.
$$
If $\pi_m(\mathbf x)=1$, then $\mathbf x=(\underbrace{0,0,0,\ldots,0}_{m-1},1,1,1,\ldots)$ and $\orbit_G(\mathbf x)=\{0,1\}$, which is not dense in $X$.  If $\pi_m(\mathbf x)=\frac{2^{n+1}-1}{2^{n+1}}$ for some non-negative integer $n$, then we distinguish the following possible cases.
\begin{enumerate}
\item $n=0$. Then $\mathbf x=(\underbrace{0,0,0,\ldots,0}_{m-1},\frac{1}{2},\frac{1}{2},\frac{1}{2},\ldots)$ and $\orbit_G(\mathbf x)=\{0,\frac{1}{2}\}$, which is not dense in $X$. 
\item $n\geq 1$. Then there is $k\in \{1,2,3,\ldots, n\}$ such that 
$$
\mathbf x=(\underbrace{0,0,0,\ldots,0}_{m-1},\frac{2^{n}-1}{2^{n}},a_k^{n},a_k^{n},a_k^{n},\ldots)
$$
 and $\orbit_G(\mathbf x)=\{0,\frac{2^{n}-1}{2^{n}},a_k^{n}\}$, which is not dense in $X$. 
\end{enumerate}
We have just proved that $0\not \in \trans_2(G)$.  Note that for any $x\in X$, there is a point $\mathbf x\in \tplus_G(0)$ such that $x\in \orbit_G(\mathbf x)$. Therefore, $\mathcal U_G^{\oplus}(0)$ is dense in $X$. It follows that $0\in \trans_3(G)$. 
Note that $\trans_3(G)=\{0\}$ and, therefore, $\trans_3(G)$ is not dense in $X$.  Also, note that $\illegal(G)\neq \emptyset$.
\end{example}
\section{3-transitive points that are not 2-transitive} \label{s2}
In this section, we study true 3-transitive points; i.e., the 3-transitive points that are not 2-transitive.  We begin with the following definition. 
\begin{definition}
Let $(X,G)$ be a CR-dynamical system and let $x\in X$.  For each positive integer $n$, we define 
$$
\mathcal R_n(x)=\{x\}\cup G(x)\cup G^{2}(x)\cup G^{3}(x)\cup \ldots \cup G^n(x)
$$
to be \emph{\color{blue} the $n$-reach of the point $x$}.  We also define
$$
\mathcal R_{\omega}(x)=\{x\}\cup G(x)\cup G^{2}(x)\cup G^{3}(x)\cup \ldots
$$
to  be \emph{\color{blue} the $\omega$-reach of the point $x$}.  
\end{definition}
\begin{observation}
Let $(X,G)$ be a CR-dynamical system and let $x\in \trans_3(G)$.  Then 
$$
\mathcal R_{\omega}(x)=\orbitt_G(x).
$$
Also, note that $x\in \trans_3(G)$ if and only if $\mathcal R_{\omega}(x)$ is dense in $X$.
\end{observation}
\begin{definition}\label{666}
Let $(X,G)$ be a CR-dynamical system and let $x\in \trans_3(G)\setminus \trans_2(G)$.  We say that $x$ is 
\begin{enumerate}
\item \emph{\color{blue}  $(3,n)$-transitive in $(X,G)$}, if there is a positive integer $m$ such that $\mathcal R_m(x)$ is dense in $X$ and 
$$
n=\min\{m \ | \ m \textup{ is a positive integer such that } \mathcal R_m(x) \textup{ is dense in } X\}.
$$
  We use \emph{\color{blue}  $\trans_{(3,n)}(G)$} to denote the set of  $(3,n)$-transitive points in $(X,G)$.
  \item \emph{\color{blue} $(3,\omega)$-transitive in $(X,G)$}, if 
  $$
  x\not \in \bigcup_{n=1}^{\infty}\trans_{(3,n)}(G).
  $$
  We use \emph{\color{blue}  $\trans_{(3,\omega)}(G)$} to denote the set of  $(3,\omega)$-transitive points in $(X,G)$.
\end{enumerate} 
\end{definition}
\begin{observation}
Let $(X,G)$ be a CR-dynamical system and let $x\in X$.  For all positive integers $m,n$,
$$
x\in \trans_{(3,n)}(G) \textup{ and } m\geq n \Longrightarrow  \mathcal R_m(x) \textup{ is dense in } X. 
$$
Also, if $x$ is illegal, then there is a positive integer $n$ such that for all integers $m\geq n$, $\mathcal R_m(x) =\mathcal R_n(x)$.
\end{observation}

\begin{example}\label{ex3a}
Let $X=[0,1]$ and let $G=([0,1]\times \{1\})\cup(\{0\}\times [0,1])$. 
We have seen in Example \ref{ex3} that $0\in \trans_3(G)$ and that for each $y\in (0,1]$, $(0,y)\in G$. It follows that 
$$
0\in \trans_{(3,1)}(G).
$$
\end{example}
Note that there are  CR-dynamical systems $(X,G)$ such that $\trans_{(3,n)}(G)\neq\emptyset$ and $\trans_{(3,n)}(G)$ is not dense in $X$; Example \ref{ex3a} is such an example.  
%
We continue by examining even more such examples.
\begin{example}\label{ex31a}
Let $X=[0,1]$, and let $G$ be defined as in Example \ref{ex31}.
We have seen in Example \ref{ex31} that $0\in \trans_3(G)\setminus \trans_2(G)$. Note that for each positive integer $n$, 
$$
0\not\in \trans_{(3,n)}(G).
$$
Therefore, 
$$
0\in \trans_{(3,\omega)}(G).
$$
Note that the  union of two orbit sets, $\orbit_G(\mathbf x_1)$ and $\orbit_G(\mathbf x_2)$ (see Example \ref{ex31} for more details) guaranties the density of the set $\orbitt_G(0)$.
\end{example}
\begin{example}\label{ex32a}
Let $(X,G)$ be the CR-dynamical system from Example \ref{ex32}.
We have seen in that example that  $0\in \trans_3(G)\setminus \trans_2(G)$. Note that for each positive integer $n$, 
$$
0\not\in \trans_{(3,n)}(G).
$$
Therefore, 
$$
0\in \trans_{(3,\omega)}(G).
$$
Note that the  union of (countably) infinitely many orbit sets guaranties the density if the set $\orbitt_G(0)$ while it is not possible to achieve this with any union of just finitely many orbit sets, see Example \ref{ex32} for more details.
\end{example}
Examples \ref{ex31a} and \ref{ex32a} also serve as a motivation for the following definition.
\begin{definition}
Let $(X,G)$ be a CR-dynamical system, let $x\in \trans_{(3,\omega)}(G)$ and let $n$ be a positive integer. We say that $x$ is
\begin{enumerate}
\item \emph{\color{blue} $(3,\omega,n)$-transitive in $(X,G)$}, if $x\in \trans_{(3,\omega)}(G)$ and there are 
$$
\mathbf x_1,\mathbf x_2,\mathbf x_3,\ldots ,\mathbf x_n\in \tplus_G(x)
$$
 such that
\begin{enumerate}
\item for all $i,j\in \{1,2,3,\ldots, n\}$, 
$$
\mathbf x_i=\mathbf x_j \Longleftrightarrow i=j,
$$
\item $\orbit_G(\mathbf x_1)\cup\orbit_G(\mathbf x_2)\cup\orbit_G(\mathbf x_3)\cup\ldots \orbit_G(\mathbf x_n)$ is dense in $X$,
\end{enumerate}
and for each positive integer $m$, if   $\mathbf x_1,\mathbf x_2,\mathbf x_3,\ldots ,\mathbf x_m\in \tplus_G(x)$ are such that
\begin{enumerate}
\item for all $i,j\in \{1,2,3,\ldots, m\}$, 
$$
\mathbf x_i=\mathbf x_j \Longleftrightarrow i=j,
$$
\item $\orbit_G(\mathbf x_1)\cup\orbit_G(\mathbf x_2)\cup\orbit_G(\mathbf x_3)\cup\ldots \orbit_G(\mathbf x_m)$ is dense in $X$,
\end{enumerate}
then $m\geq n$.
  We use \emph{\color{blue}  $\trans_{(3,\omega,n)}(G)$} to denote the set of  $(3,\omega,n)$-transitive points in $(X,G)$.
  \item \emph{\color{blue} $(3,\omega,\omega)$-transitive in $(X,G)$}, if $x\in \trans_{(3,\omega)}(G)$ and for each positive integer $n$, $x\not \in \trans_{(3,\omega ,n)}(G)$. 
  We use \emph{\color{blue}  $\trans_{(3,\omega,\omega)}(G)$} to denote the set of  $(3,\omega,\omega)$-transitive points in $(X,G)$.
\end{enumerate}
\end{definition}
We conclude this section with the following two examples.
\begin{example}\label{ex31aa}
Let $(X,G)$ be the CR-dynamical system from Example \ref{ex31}.
Then $0\in \trans_{(3,\omega,2)}(G)$.
\end{example}
\begin{example}\label{ex32aa}
Let $(X,G)$ be the CR-dynamical system from Example \ref{ex32}.
Then $0\in \trans_{(3,\omega,\omega)}(G)$.
\end{example}
\section{Transitivity trees}  \label{s3}
In the previous section, we saw that 3-transitivity leads to many different possible and quite complicated orbit structures.  These structures naturally form trees of countable height.  

In this section,  we present the notion of transitive points by a set-theoretical model, a tree, as a tool to help provide an intuitive (visual) guide to transitivity in this setting, to potentially aid future research.

We restrict our definition of trees to connected trees of height $\omega$. We employ $\mathbb N$ in place of the ordinal $\omega$.
\begin{definition}
Let $(T,\leq)$ be a partially ordered set.  We say that $(T,\leq)$ is a \emph{\color{blue}  tree}, if there is a unique point $x\in T$, the root of $T$,  such that 
\begin{enumerate}
\item for each $y\in T$, $x\leq y$ and 
\item for each $y\in T$,   $(\{z\in T \ | \ z\leq y\},\leq)$ is well ordered. 
\end{enumerate}
Maximal well-ordered sets in a tree are called \emph{\color{blue} branches}.  If $n$ is the cardinality of a branch $B$, then we say that \emph{\color{blue} the length of the branch $B$ is $n-1$} and write \emph{\color{blue} $\ell_T(B)=n-1$}. If $n$ is finite, we say that \emph{\color{blue} the length of $B$ is finite} and write \emph{\color{blue}$\ell_T(B)<\infty$}.  If $n$ is not finite, we say that \emph{\color{blue} the length of $B$ is infinite} and write \emph{\color{blue}$\ell_T(B)=\infty$}.  We denote by  \emph{\color{blue} $\mathcal B(T)$} the set of branches of the tree $(T,\leq)$ and by  \emph{\color{blue} $\mathcal B_{\infty}(T)$} the set of branches of the tree $(T,\leq)$ with infinite length. We also define  \emph{\color{blue}  the height of $(T,\leq)$ } as
$$
\textup{height}(T)=\sup\{\ell_T(B) \ | \ B \textup{ is a branch of } (T,\leq)\}.
$$
\end{definition}

\begin{figure}[h!]
	\centering
		\includegraphics[width=21em]{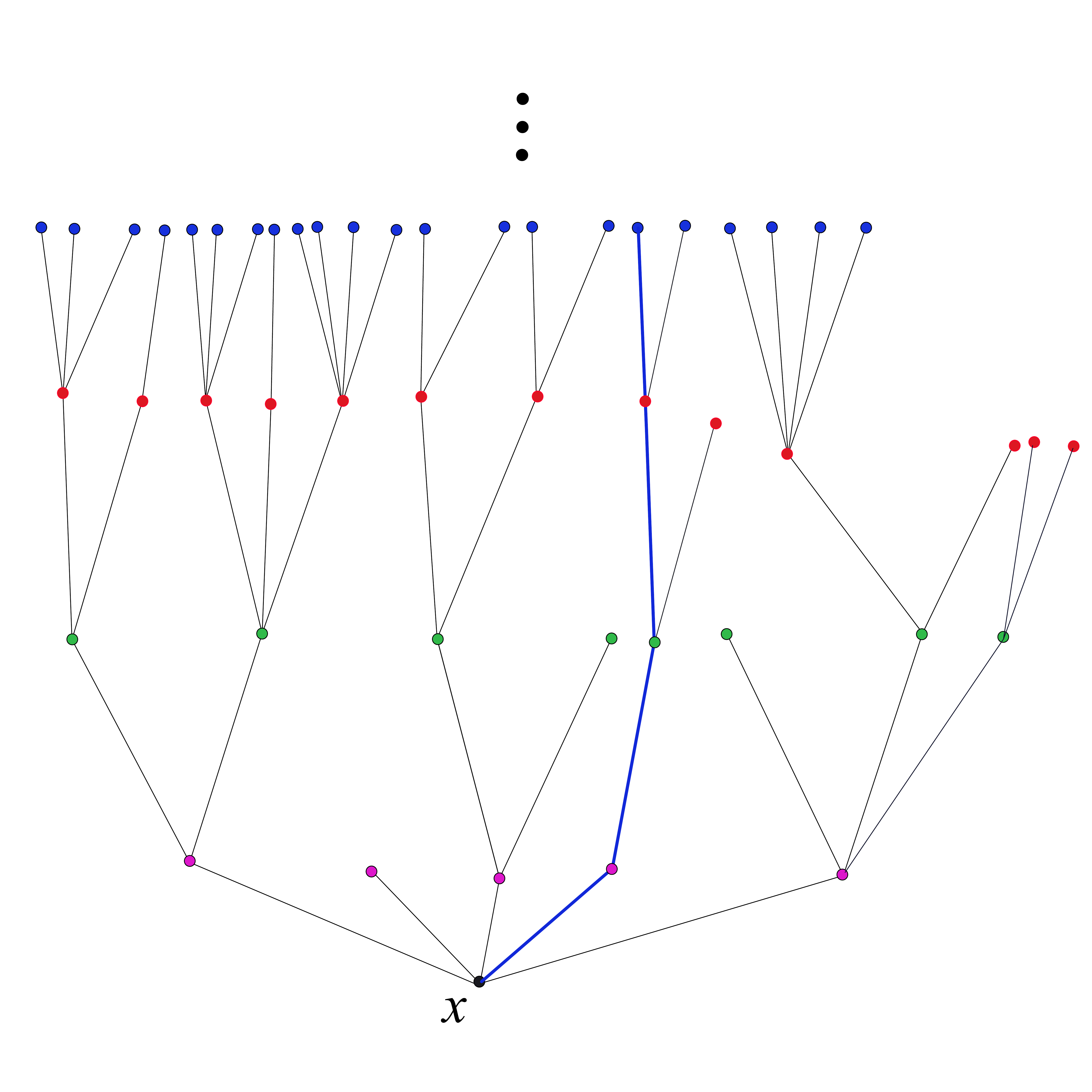}
	\caption{A tree}
	\label{tree}
\end{figure}

\begin{definition}
Let $(T,\leq)$ be a tree with root $x\in T$. We define  \emph{\color{blue} $\textup{level}_0(T)$} by
$$
\textup{level}_0(T)=\{x\}.
$$
Let $n$ be a positive integer.  Then we say that  \emph{\color{blue} $\textup{level}_n(T)$} is the set of all points $y\in T$ such that there are points $z_0,z_1,z_2,z_3,\ldots ,z_n\in T$ satisfying 
\begin{enumerate}
\item $z_0=x$ and $z_n=y$,
\item for each $i\in\{0,1,2,3,\ldots, n-1\}$, 
$$
z_i\leq z_{i+1} \textup{ and } z_i\neq z_{i+1}
$$
\item for each $i\in\{0,1,2,3,\ldots, n-1\}$ and for each $w\in T$, 
$$
z_i\leq w\leq z_{i+1} \Longrightarrow  w\in \{z_i, z_{i+1}\}.
$$
\end{enumerate}  
\end{definition}
In Figure \ref{tree}, there is an example of a tree $(T,\leq)$ with root $x$.  A branch of the tree is coloured blue.  The points in $\level_1(T)$ are coloured pink, the points in $\level_2(T)$ are coloured green, the points in $\level_3(T)$ are coloured red, and the points in $\level_4(T)$ are coloured blue.  
\begin{definition}
Let $(T,\leq)$ be a tree with  root $x\in T$ and let $n$ be a non-negative integer. We define the set  \emph{\color{blue}$\mathcal L_n(T)$} as follows:
$$
\mathcal L_n(T)=\level_0(T)\cup \level_1(T)\cup \level_2(T)\cup\ldots \cup \level_n(T).
$$
\end{definition}
Next, we implement the trees into the theory of CR-dynamical systems.
\begin{definition}
Let $(X,G)$ be a CR-dynamical system and let $x\in X$. Then we define  \emph{\color{blue} $(T(x),\leq)$} as follows.  Let 
$$
T(x)=\{x\}\cup\bigcup_{n=1}^{\infty}G^n(x)
$$
and for all $y,z\in X$, we define
$$
y\leq z \Longleftrightarrow z=y  \textup{ or  there is a positive integer } n \textup{ such that } z\in G^n(y).
$$
\end{definition}
\begin{observation}
Let $(X,G)$ be a CR-dynamical system and let $x\in X$. Then $(T(x),\leq)$ is a tree.
\end{observation}
\begin{definition}
Let $(X,G)$ be a CR-dynamical system and let $x\in X$. We call the tree $(T(x),\leq)$ \emph{\color{blue} the transitivity tree of $(X,G)$ with respect to $x$}.
\end{definition}
\begin{observation}
Let $(X,G)$ be a CR-dynamical system, let $x\in X$, and let $(T(x),\leq)$ be the transitivity tree of $(X,G)$ with respect to $x$.  The leaves of the tree $(T(x),\leq)$ are the maximal members of branches and they correspond to illegal points.  

Also,  the point $x$ is illegal if and only if $(T(x),\leq)$ is a tree of finite height. 
\end{observation}
\begin{observation}
Let $(X,G)$ be a CR-dynamical system, let $x\in X$, let $n$ be a positive integer,  and let $(T(x),\leq)$ be the transitivity tree of $(X,G)$ with respect to $x$.  Then the following holds.
\begin{enumerate}
\item $x\in \legal(G)$ if and only if $\mathcal B_{\infty}(T(x))\neq \emptyset$.
\item $x\in trans_1(G)$ if and only if  $\mathcal B_{\infty}(T(x))\neq \emptyset$ and for each $B\in \mathcal B_{\infty}(T(x))$,  $\Cl(B)=X$.
\item $x\in trans_2(G)$ if and only if there is  $B\in \mathcal B_{\infty}(T(x))$ such that 
$\Cl(B)=X$.
\item $x\in trans_3(G)$ if and only if $ \Cl\left(\bigcup_{B\in \mathcal B_{\infty}(T(x))}B\right)=X$.
\item $x\in \intrans(G)$ if and only if $\mathcal B_{\infty}(T(x))\neq \emptyset$ and $ \Cl\left(\bigcup_{B\in \mathcal B_{\infty}(T(x))}B\right)\neq X$.
\item  $x\in \trans_{(3,n)}(G)$ if and only if $\Cl(\mathcal L_n(T(x)))=X$ and $\Cl(\mathcal L_{n-1}(T(x)))\neq X$.
\item \label{cmuz1} $x\in \trans_{(3,\omega)}(G)$ if and only if $\Cl(T(x))=X$, for each $B\in \mathcal B_{\infty}(T(x))$, $\Cl(B)\neq X$,  and for each positive integer $n$, $\Cl(\mathcal L_n(T(x)))\neq X$.  
\item\label{cmuz2}  $x\in \trans_{(3,\omega,n)}(G)$ if and only if  $x\in \trans_{(3,\omega)}(G)$ and there are $T_1,T_2,T_3,\ldots ,T_n\in \mathcal B_{\infty}(T(x))$  such that
\begin{enumerate}
\item for all $i,j\in \{1,2,3,\ldots, n\}$, 
$$
T_i=T_j \Longleftrightarrow i=j,
$$
\item $T_1\cup T_2\cup T_3\cup \ldots \cup T_n$ is dense in $X$,
\end{enumerate}
and for each positive integer $m$, if $T_1,T_2,T_3,\ldots ,T_m\in \mathcal B_{\infty}(T(x))$  such that
\begin{enumerate}
\item for all $i,j\in \{1,2,3,\ldots, m\}$, 
$$
T_i=T_j \Longleftrightarrow i=j,
$$
\item $T_1\cup T_2\cup T_3\cup \ldots \cup T_m$ is dense in $X$,
\end{enumerate}
then $m\geq n$.
\item $x\in \trans_{(3,\omega,\omega)}(G)$ if and only if statement \ref{cmuz1} is satisfied and statement \ref{cmuz2} is not.
\end{enumerate}
\end{observation}
\begin{observation}
Let $(X,G)$ be a CR-dynamical system. Then the following hold.
\begin{enumerate}
\item $G$ is a graph of a single valued function from a subspace of $X$ to $X$ if and only if for each $x\in X$,  $|\mathcal B(T(x))|=1$.
\item $G$ is a graph of a single valued function from $X$ to $X$ if and only if for each $x\in X$,   $|\mathcal B_{\infty}(T(x))|=1$.
\end{enumerate}
\end{observation}
\section{Dense orbit transitive CR-dynamical systems} \label{s4}
We begin the final section by recalling the definition of a dense orbit transitive dynamical system $(X,f)$. 
\begin{definition}
Let $(X,f)$ be a dynamical system.  We say that $(X,f)$ is 
 \emph{ \color{blue} dense orbit transitive} or \emph{ \color{blue} DO-transitive},   if $\tr(f)\neq \emptyset$.
\end{definition}
Next, we generalize this definition to CR-dynamical systems.  Since in CR-dynamical systems, there are many different types of transitive points, it is natural to introduce different types of dense orbit transitivity in this context, see Definition \ref{DEDO}.
\begin{definition}\label{DEDO}
Let $(X,G)$ be a CR-dynamical system.  We say that 
\begin{enumerate}
\item  $(X,G)$ is \emph{\color{blue} type 1 dense orbit transitive} or \emph{\color{blue} 1-DO-transitive}, if $\trans_1(G)\neq \emptyset$.
\item  $(X,G)$ is \emph{\color{blue} type 2 dense orbit transitive} or \emph{\color{blue} 2-DO-transitive}, if $\trans_2(G)\neq \emptyset$.
\item  $(X,G)$ is \emph{\color{blue} type 3 dense orbit transitive} or \emph{\color{blue} 3-DO-transitive}, if $\trans_3(G)\neq \emptyset$.
\item  For each positive integer $n$, we say that $(X,G)$ is \emph{\color{blue} type $(3,n)$ dense orbit transitive} or \emph{\color{blue} $(3,n)$-DO-transitive}, if $\trans_{(3,n)}(G)\neq \emptyset$. 
\item We say that $(X,G)$ is \emph{\color{blue} type $(3,\omega)$ dense orbit transitive} or \emph{\color{blue} $(3,\omega)$-DO-transitive}, if $\trans_{(3,\omega)}(G)\neq \emptyset$.
\item  For each positive integer $n$, we say that $(X,G)$ is \emph{\color{blue} type $(3,\omega,n)$ dense orbit transitive} or \emph{\color{blue} $(3,\omega,n)$-DO-transitive}, if $\trans_{(3,\omega,n)}(G)\neq \emptyset$. 
\item We say that $(X,G)$ is \emph{\color{blue} type $(3,\omega,\omega)$ dense orbit transitive} or \emph{\color{blue} $(3,\omega,\omega)$-DO-transitive}, if $\trans_{(3,\omega,\omega)}(G)\neq \emptyset$.
\end{enumerate}
\end{definition}
\begin{observation}\label{observation1}
Let $(X,G)$ be a CR-dynamical system. If $(X,G)$ is 1-DO-transitive, then $(X,G)$ is 2-DO-transitive.  If $(X,G)$ is 2-DO-transitive, then $(X,G)$ is 3-DO-transitive. 
\end{observation}
\begin{observation}
Let $(X,f)$ be a dynamical system. The following statements are equivalent.
\begin{enumerate}
\item $(X,f)$ is DO-transitive.
\item $(X,\Gamma(f))$ is 1-DO-transitive.
\item $(X,\Gamma(f))$ is 2-DO-transitive.
\item $(X,\Gamma(f))$ is 3-DO-transitive.
\end{enumerate}
\end{observation}
It is a well-known fact that if $(X,f)$ is a DO-transitive dynamical system and  $\isolated(X)=\emptyset$, then $f$ is surjective.  Theorem \ref{erceg} generalizes this to CR-dynamical systems. Also, note that there are examples of DO-transitive dynamical system (with  $\isolated(X)\neq \emptyset$) such that $f$ is not surjective; for an example see Example \ref{tistile} letting $f:X\rightarrow X$ be such a mapping that $G=\Gamma(f)$. 
\begin{theorem}\label{erceg}
Let $(X,G)$ be a CR-dynamical system,  such that $X$ has no isolated points or $X$ is degenerate.  Then for each $k\in \{1,2,3\}$,
$$
(X,G) \textup{ is k-DO-transitive } \Longrightarrow p_1(G)=p_2(G)=X.
$$
\end{theorem}
\begin{proof}
If $X$ is degenerate, then there is nothing to prove.  For the rest of the proof, assume that $X$ is a non-degenerate space.
Let $(X,G)$ be a 3-DO-transitive CR-dynamical system and let $x\in \trans_3(G)$. Then $\orbitt_G(x)$ is dense in $X$.  To see that $p_1(G)=p_2(G)=X$, let $t\in X$ be any point.  We show that $t\in p_1(G)\cap p_2(G)$.  First, we prove the following claim.
\begin{nonumclaim}
There is a sequence  of points $(t_n)$ in $\orbitt_G(x)$ such that
\begin{enumerate}
\item $\displaystyle \lim_{n\to \infty}t_n=t$ and 
\item for each positive integer $n$, there is a positive integer $j_n\geq n$ such that $t_{j_n}\neq t$.
\end{enumerate} 
\end{nonumclaim}
\begin{proof}[Proof of Claim]\renewcommand{\qed}{} 
Suppose that for each sequence $(t_n)$ of points in $\orbitt_G(x)$ such that $\displaystyle \lim_{n\to \infty}t_n=t$, there is a positive integer $n_0$ such that for each $n\geq n_0$, $t_n=t$.  First, we show that there is an open set $U$ in $X$ such that 
$$
U\cap \orbitt_G(x)=\{t\}.
$$ 
On the contrary, suppose that for each open set $U$ in $X$ such that $t\in U$,  there is $z\in  U\cap \orbitt_G(x)$ such that $z\neq t$.  Using that, let for each positive integer $n$,  
$$
t_n\in \Big(B\Big(t,\frac{1}{n}\Big)\cap \orbitt_G(x)\Big)\setminus\{t\}. 
$$
Then $(t_n)$ is a sequence in $\orbitt_G(x)$ such that $\displaystyle \lim_{n\to \infty}t_n=t$ for each positive integer $n$,  $t_{n}\neq t$ -- a contradiction.  Therefore,  there is an open set $U$ in $X$ such that 
$$
U\cap \orbitt_G(x)=\{t\}.
$$ 
Let $U$ be such an open set in $X$. Since $X$ has no isolated points, $U\neq \{t\}$.  
Then $U\setminus \{t\}$ is a non-empty open set in $X$. Since $\orbitt_G(x)$ is dense in $X$, it follows that $(U\setminus \{t\})\cap \orbitt_G(x)\neq \emptyset$ -- a contradiction.
This completes the proof of the claim.
\end{proof}
It follows from the above claim that there is a sequence  of points $(t_n)$ in $\orbitt_G(x)$ such that
\begin{enumerate}
\item $\displaystyle \lim_{n\to \infty}t_n=t$ and 
\item for each positive integer $n$,  $t_{n}\neq t$.
\end{enumerate} 
Let $(t_n)$ be such a sequence.  For each positive integer $n$, let $\mathbf x_n\in \tplus_G(x)$ and let $i_n$ be a positive integer such that $\pi_{i_n}(\mathbf x_n)=t_{n}$. We consider the following possible cases.
\begin{enumerate}
\item[(i)] $x\neq t$. Then there is a positive integer $n_0$ such that for each positive integer $n\geq n_0$, $t_{n}\neq x$. It follows that for each positive integer $n\geq n_0$, $i_n>1$. It follows that for each positive integer $n\geq n_0$, 
$$
(\pi_{i_n-1}(\mathbf x_n),\pi_{i_n}(\mathbf x_n),\pi_{i_n+1}(\mathbf x_n))\in \star_{i=1}^{2}G.
$$
 Note that for each positive integer $n$, $\pi_{i_n}(\mathbf x_n)=t_{n}$ and that $\star_{i=1}^{2}G$ is compact.  Let $(x_0,t_0,y_0)\in \star_{i=1}^{2}G$ and let $(\pi_{i_{k_n}-1}(\mathbf x_{k_n}),t_{{k_n}},\pi_{i_{k_n}+1}(\mathbf x_{k_n}))$ be a convergent subsequence of the sequence $(\pi_{i_n-1}(\mathbf x_n),t_{n},\pi_{i_n+1}(\mathbf x_n))$ such that
$$
\lim_{n\to \infty}(\pi_{i_{k_n}-1}(\mathbf x_{k_n}),t_{{k_n}},\pi_{i_{k_n}+1}(\mathbf x_{k_n}))=(x_0,t_0,y_0).
$$
Then $t_0=t$ and it follows that $t\in p_1(G)\cap p_2(G)$. 
\item[(ii)] $x=t$. 
Since for each positive integer $n$, $\pi_{i_n}(\mathbf x_n)=t_{n}$, it follows that for each positive integer $n$, $i_n>1$. Therefore, for each positive integer $n$, 
$$
(\pi_{i_n-1}(\mathbf x_n),\pi_{i_n}(\mathbf x_n),\pi_{i_n+1}(\mathbf x_n))\in \star_{i=1}^{2}G.
$$
Following the proof of (i) for $n_0=1$, we also get here that $t\in p_1(G)\cap p_2(G)$. 
\end{enumerate}

Next, let $(X,G)$ be a 1-DO-transitive or 2-DO-transitive CR-dynamical system. By Observation \ref{observation1},  $(X,G)$ is a 3-DO-transitive CR-dynamical system. It follows that $p_1(G)=p_2(G)=X$.
\end{proof}

It is trivial to note that when dealing with a standard dynamical systems $(X,f)$ it is always the case that $p_1(G)=X$, if $G=\Gamma(f)$. However,  it is not true for any CR-dynamical system $(X,G)$ that $p_1(G)=X$.  In Theorem \ref{gerce}, we show that if $(X,G)$ is a $k$-DO-transitive CR-dynamical system ($k\in\{1,2,3\}$), then $p_1(G)=X$ no matter if $\isolated (X)=\emptyset$ or $\isolated(X)\neq \emptyset$.  
	\begin{theorem}\label{gerce}
		Let $(X,G)$ be a CR-dynamical system. Then for each $k\in \{1,2,3\}$,
		$$
		(X,G) \textup{ is k-DO-transitive } \Longrightarrow p_1(G)=X.
		$$
	\end{theorem}
	\begin{proof}
{	First, suppose that $(X,G)$ is 3-DO-transitive and let $x\in \trans_3(G)$. Then $\orbitt_G(x)$ is dense in $X$, it means that 
	$$
	\Cl(\orbitt_G(x))=X.
	$$
	  To see that $p_1(G)=X$, let $t\in X$ be any point.  We show that $t\in p_1(G)$.  Since $\Cl(\orbitt_G(x))=X$,  it follows that for each positive integer $n$, there is 
	  $$
	  t_n\in B\Big(x,\frac{1}{n}\Big)\cap \orbitt_G(x).
	  $$
	  For each positive integer $n$, fix such a point $t_n$.  Since $\displaystyle \orbitt_G(x)=\bigcup_{\mathbf y\in T_G^{+}(x)}\orbit_G(\mathbf y)$, it follows that for each positive integer $n$,  there is $\mathbf y_n\in T_G^{+}(x)$ such that
	  $$
t_n\in 	  B\Big(x,\frac{1}{n}\Big)\cap \orbit_G(\mathbf y_n).
	  $$
	  Since for each positive integer $n$, $\orbit_G(\mathbf y_n)\subseteq p_1(G)$, it follows that for each positive integer $n$, 
	  $$
	  t_n\in B\Big(x,\frac{1}{n}\Big)\cap p_1(G).
	  $$
	Since $p_1(G)$ is closed in $X$ and since $\displaystyle \lim_{n\to \infty}t_n=t$, it follows that $t\in p_1(G)$.
	
	Next, suppose that $(X,G)$ is 1-DO-transitive or 2-DO-transitive. It follows from Observation \ref{observation1} that $(X,G)$ is also 3-DO-transitive and we have just proved that in this case $p_1(G)=X$ follows. This completes the proof. }
	\end{proof}
	We conclude this section with the following example.
	\begin{example}\label{tistile}
		Let $X=[0,1]\cup\{2\}$ and let $f:[0,1]\rightarrow [0,1]$ be the tent-map defined by $f(t)=2t$ for $t\leq \frac{1}{2}$ and $f(t)=2-2t$ for $t\geq \frac{1}{2}$,  let $x\in \tr(f)$ and let  $G=\Gamma(f)\cup \{(2,x)\}$, see Figure \ref{figgure}.
		\begin{figure}[h!]
			\centering
			\includegraphics[width=20em]{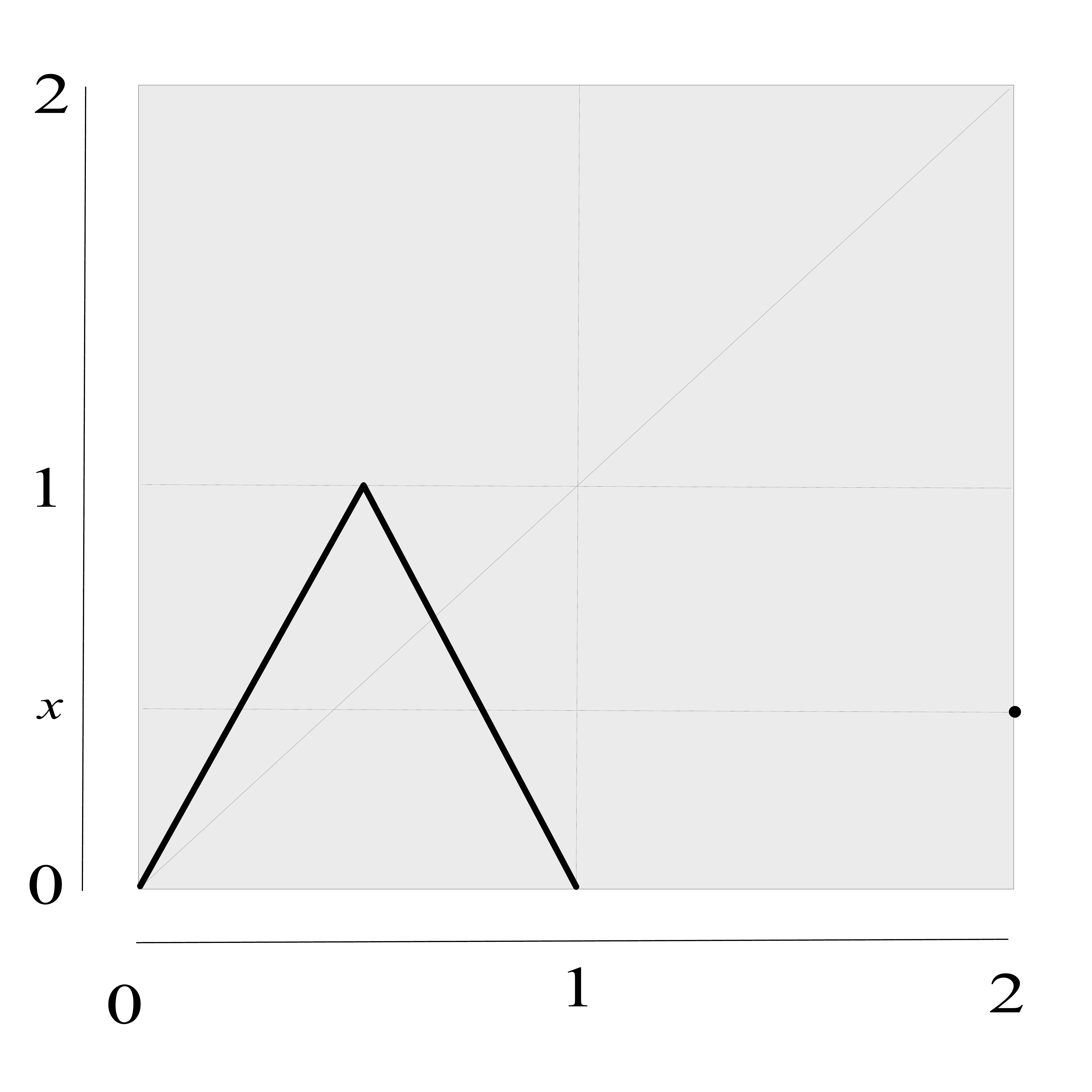}
			\caption{The relation  $G$ from Example \ref{tistile}} 
			\label{figgure}
		\end{figure} 
\noindent Note that $\trans_1(G)=\trans_2(G)=\trans_3(G)=\{2\}$ and that $p_2(G)\neq X$.
	\end{example}

\section{Transitive CR-dynamical systems}\label{s5}
First, we revisit the definition of a transitive dynamical system. 
\begin{definition}\label{simpl}
Let $(X,f)$ be a dynamical system.  We say that $(X,f)$ is 
 \emph{\color{blue} transitive}, if for all non-empty open sets $U$ and $V$ in $X$,  there is a non-negative integer $n$ such that 
$$
f^n(U)\cap V\neq \emptyset.
$$
\end{definition}
Theorem \ref{jurc} is a well known result, see \cite{B,KS} for more information.
\begin{theorem}\label{jurc}
Let $(X,f)$ be a dynamical system. The following statements are equivalent.
\begin{enumerate}
\item$(X,f)$ is  transitive.
\item For all non-empty open sets $U$ and $V$ in $X$,  there is a positive integer $n$ such that 
$$
f^{n}(U)\cap V\neq \emptyset.
$$
\item  For each non-empty open set $U$ in $X$,  $\bigcup_{k=0}^{\infty}f^{k}(U)$ is dense in $X$.
\item  For each non-empty open set $U$ in $X$,  $\bigcup_{k=1}^{\infty}f^{k}(U)$ is dense in $X$.
\item  For all non-empty open sets $U$ and $V$ in $X$,  there is a non-negative integer $n$ such that 
$$
f^{-n}(U)\cap V\neq \emptyset.
$$
\item For all non-empty open sets $U$ and $V$ in $X$,  there is a positive integer $n$ such that 
$$
f^{-n}(U)\cap V\neq \emptyset.
$$
\item  For each non-empty open set $U$ in $X$,  $\bigcup_{k=0}^{\infty}f^{-k}(U)$ is dense in $X$.
\item  For each non-empty open set $U$ in $X$,  $\bigcup_{k=1}^{\infty}f^{-k}(U)$ is dense in $X$.
\end{enumerate}
\end{theorem}
In the following definition, we  generalize  Definition \ref{simpl} and in Theorem \ref{transitiveCHARACTERIZATION11}, we generalize Theorem \ref{jurc} to CR-dynamical systems.
\begin{definition}
Let $(X,G)$ be a CR-dynamical system.  We say that $(X,G)$ is 
 \emph{\color{blue} transitive}, if for all non-empty open sets $U$ and $V$ in $X$,  there is a non-negative integer $n$ such that 
$$
G^n(U)\cap V\neq \emptyset.
$$
\end{definition}
We use the following lemma in the proof of Theorem \ref{transitiveCHARACTERIZATION11}.
\begin{lemma}\label{lemai}
Let $(X,G)$ be a CR-dynamical system. If $X$ has no isolated points, then \ref{on1} is equivalent to  \ref{tw2}.
\begin{enumerate}
\item\label{on1} $(X,G)$ is transitive.
\item\label{tw2} For any non-empty  open set $U$ in $X$, 
$$
\Cl\left(\bigcup_{k=1}^{\infty}G^{-k}(U)\right)=X.
$$
\end{enumerate}
\end{lemma}
\begin{proof}
Suppose that $(X,G)$ is transitive. 
Let $U$ be a non-empty open subset of $X$. To show that $\Cl\left(\bigcup_{k=1}^{\infty}G^{-k}(U)\right)=X$, let $V$ be a non-empty open subset of $X$.  We show that 
$$
V\cap \left(\bigcup_{k=1}^{\infty}G^{-k}(U)\right)\neq \emptyset.
$$
We treat the following possible cases.
\begin{enumerate}
\item $V\setminus \Cl(U)\neq \emptyset$.  Let $W=V\setminus \Cl(U)$. Then $W$ is a non-empty open subset of $X$.  Next, let $m$ be a non-negative integer such that 
$$
G^m(W)\cap U\neq \emptyset,
$$
and let $\mathbf x=(x_1,x_2,x_3,\ldots, x_m,x_{m+1})\in \star_{i=1}^{m}G$ such that $x_1\in W$ and $x_{m+1}\in U$.  Note that $m\geq 1$.  Then 
$$
(x_{m+1},x_m,x_{m-1},\ldots, x_2,x_{1})\in \star_{i=1}^{m}G^{-1}
$$
and it follows that $x_{1}\in V\cap G^{-m}(U)$. Therefore, 
$
V\cap \left(\bigcup_{k=1}^{\infty}G^{-k}(U)\right)\neq \emptyset.
$
 \item $V\setminus \Cl(U)=\emptyset$.   It follows that $V\subseteq \Cl(U)$.  Since $X$ has no isolated points, it follows that $V$ has at least two points.  Let $x,y\in V$ such that $x\neq y$.  Let 
 $$
 W_1=B\Big(x,\frac{d(x,y)}{2}\Big)\cap V \textup{ and } W_2=B\Big(y,\frac{d(x,y)}{2}\Big)\cap V.
 $$
 Then $W_1$ and $W_2$ are non-empty open subsets of $X$ such that 
 $$
 W_1\cap W_2=\emptyset \textup{ and } W_1\cup W_2\subseteq V.
 $$
 It follows from $W_1\subseteq \Cl(U)$ and $W_2\subseteq \Cl(U)$ that $W_1\cap U\neq \emptyset$ and $W_2\cap U\neq \emptyset$. Let $m$ be a non-negative integer such that 
 $$
 G^m(W_1\cap U)\cap(W_2\cap U)\neq \emptyset.
 $$ 
 Such a non-negative integer $m$ does exist since $(X,G)$ is transitive.  Let $\mathbf x=(x_1,x_2,x_3,\ldots, x_m,x_{m+1})\in \star_{i=1}^{m}G$ such that $x_1\in W_1\cap U$ and $x_{m+1}\in W_2\cap U$.  Note that $m\geq 1$ since $W_1\cap W_2= \emptyset$.   It follows that
 $$
 x_1\in V\cap G^{-m}(x_{m+1})\subseteq V\cap G^{-m}(U)\subseteq V\cap \left(\bigcup_{k=1}^{\infty}G^{-k}(U)\right)
 $$
 and therefore $V\cap \left(\bigcup_{k=1}^{\infty}G^{-k}(U)\right)\neq \emptyset$.
\end{enumerate}
This proves the implication from \ref{on1} to \ref{tw2}. To prove the implication  from \ref{tw2} to \ref{on1}, let $U$ and $V$ be non-empty open sets in $X$.  Then $\Big(\bigcup_{k=1}^{\infty}G^{-k}(V)\Big)\cap U\neq \emptyset$.  Let $x\in  \Big(\bigcup_{k=1}^{\infty}G^{-k}(V)\Big)\cap U$ and let $n$ be a positive integer such that $x\in G^{-n}(V)\cap U$. Next, let
$$
\mathbf x=(x_1,x_2,x_3,\ldots ,x_n,x_{n+1})\in \star_{i=1}^{n}G
$$
such that $x_{n+1}\in V$ and $x_1=x$.  Since $x_1\in U$, it follows that $x_{n+1}\in G^n(U)$. Therefore, $G^{n}(U)\cap V\neq \emptyset$ and this completes the proof.
\end{proof}

\begin{theorem}\label{transitiveCHARACTERIZATION11}
Let $(X,G)$ be a CR-dynamical system.  Consider the following statements. 
\begin{enumerate}
\item \label{111} $(X,G)$ is  transitive.
\item \label{122} For all non-empty open sets $U$ and $V$ in $X$,  there is a positive integer $n$ such that 
$$
G^{n}(U)\cap V\neq \emptyset.
$$
\item \label{133} For each non-empty open set $U$ in $X$,  $\bigcup_{k=0}^{\infty}G^{k}(U)$ is dense in $X$.
\item \label{144} For each non-empty open set $U$ in $X$,  $\bigcup_{k=1}^{\infty}G^{k}(U)$ is dense in $X$.
\item \label{155} For all non-empty open sets $U$ and $V$ in $X$,  there is a non-negative integer $n$ such that 
$$
G^{-n}(U)\cap V\neq \emptyset.
$$
\item \label{166} For all non-empty open sets $U$ and $V$ in $X$,  there is a positive integer $n$ such that 
$$
G^{-n}(U)\cap V\neq \emptyset.
$$
\item \label{177} For each non-empty open set $U$ in $X$,  $\bigcup_{k=0}^{\infty}G^{-k}(U)$ is dense in $X$.
\item \label{188} For each non-empty open set $U$ in $X$,  $\bigcup_{k=1}^{\infty}G^{-k}(U)$ is dense in $X$.
\end{enumerate}
Then the following holds.
\begin{itemize}
\item The statements \ref{111}, \ref{133}, \ref{155} and \ref{177} are equivalent.
\item The statements \ref{122}, \ref{144}, \ref{166} and \ref{188} are equivalent.
\item If $X$ has no isolated points, then all statements are equivalent.
\end{itemize}
\end{theorem}

\begin{proof}
First, we prove the implication from \ref{188} to \ref{166}. Let $U$ and $V$ be non-empty open subsets in $X$. Since $\bigcup_{k=1}^{\infty}G^{-k}(U)$ is dense in $X$, 
$$
\bigcup_{k=1}^{\infty}G^{-k}(U)\cap V\neq \emptyset. 
$$ 
Then \ref{166} follows. The prove of the implication from \ref{177} to \ref{155} is analogous.

Next, we prove the implication from \ref{166} to \ref{144}. Let $U$ be a non-empty open subset in $X$. To show that $\bigcup_{k=1}^{\infty}G^{k}(U)$ is dense in $X$, let $V$ be a non-empty open set in $X$ and let $n$ be a positive integer such that $G^{-n}(V)\cap U\neq\emptyset$. Also, let $x\in G^{-n}(V)\cap U$ and let
$$
\mathbf x=(x_1,x_2,x_3,\ldots ,x_n,x_{n+1})\in \star_{i=1}^{n}G
$$
such that $x_{n+1}\in V$ and $x_1=x$.  Since $x_1\in U$, it follows that $x_{n+1}\in G^n(U)$. Therefore, $\bigcup_{k=1}^{\infty}G^{k}(U)\cap V\neq \emptyset$ 
and \ref{144} follows. The prove of the implication from \ref{155} to \ref{133} is analogous.

Next, we prove the implication from \ref{144} to \ref{122}. Let $U$ and $V$ be non-empty open subsets in $X$. Since $\bigcup_{k=1}^{\infty}G^{k}(U)$ is dense in $X$, it follows that there is a positive integer $n$ such that $G^{n}(U)\cap V\neq \emptyset$ and \ref{122} follows. The prove of the implication from \ref{133} to \ref{111} is analogous.

Next, we prove the implication from \ref{122} to \ref{188}. Let $U$ be a non-empty open subset in $X$. To show that $\bigcup_{k=1}^{\infty}G^{-k}(U)$ is dense in $X$, let $V$ be a non-empty open set in $X$ and let $n$ be a positive integer such that $G^{n}(V)\cap U\neq\emptyset$. Also, let $x\in G^{n}(V)\cap U$ and let
$$
\mathbf x=(x_1,x_2,x_3,\ldots ,x_n,x_{n+1})\in \star_{i=1}^{n}G^{-1}
$$
such that $x_{n+1}\in V$ and $x_1=x$.  Since $x_1\in U$, it follows that $x_{n+1}\in G^{-n}(U)$. Therefore, $\bigcup_{k=1}^{\infty}G^{-k}(U)\cap V\neq \emptyset$ 
and \ref{188} follows. The prove of the implication from \ref{111} to \ref{177} is analogous.

We have just proved that the statements \ref{111}, \ref{133}, \ref{155} and \ref{177} are equivalent and that the statements \ref{122}, \ref{144}, \ref{166} and \ref{188} are equivalent. 

To conclude the proof, assume that $X$ has no isolated points.  To show that the statements \ref{111},  \ref{122},  \ref{133},   \ref{144},  \ref{155},  \ref{166},  \ref{177} and \ref{188} are equivalent, it suffices to see that \ref{111} is equivalent to \ref{188}. This follows from Lemma \ref{lemai}.
\end{proof}
{
\begin{observation}\label{naminusena}
Let $(X,G)$ be a CR-dynamical system. It follows from Theorem \ref{transitiveCHARACTERIZATION11} that $(X,G)$ is transitive if and only if $(X,G^{-1})$ is transitive.
\end{observation}
}
In the following example, we show that the statements \ref{111} -- \ref{188} from Theorem  \ref{transitiveCHARACTERIZATION11} may not be equivalent if $X$ has an isolated point. Explicitly, we show that there is a CR-dynamical system $(X,G)$, which satisfies \ref{133} but not \ref{144}.  We need the following definition.
\begin{definition}
Let $(X,f)$ be a dynamical system. We say that $f$ is locally eventually onto if for each non-empty open set $U$ in $X$ there is a positive integer $n$ such that $f^n(U)=X$.
\end{definition}
\begin{observation}
Let $X=[0,1]$ and let $f:[0,1]\rightarrow [0,1]$ be the tent-map defined by $f(t)=2t$ for $t\leq \frac{1}{2}$ and $f(t)=2-2t$ for $t\geq \frac{1}{2}$. Then $f$ is locally eventually onto.
\end{observation}
\begin{example}\label{exxi}
Let $X=[0,1]\cup\{2\}$ and let $f:[0,1]\rightarrow [0,1]$ be the tent-map defined by $f(t)=2t$ for $t\leq \frac{1}{2}$ and $f(t)=2-2t$ for $t\geq \frac{1}{2}$,  let $x\in \tr(f)$ and let  
$$
G=\Gamma(f)\cup\{(2,x),(0,2)\},
$$
see Figure \ref{fi}.
\begin{figure}[h!]
	\centering
		\includegraphics[width=20em]{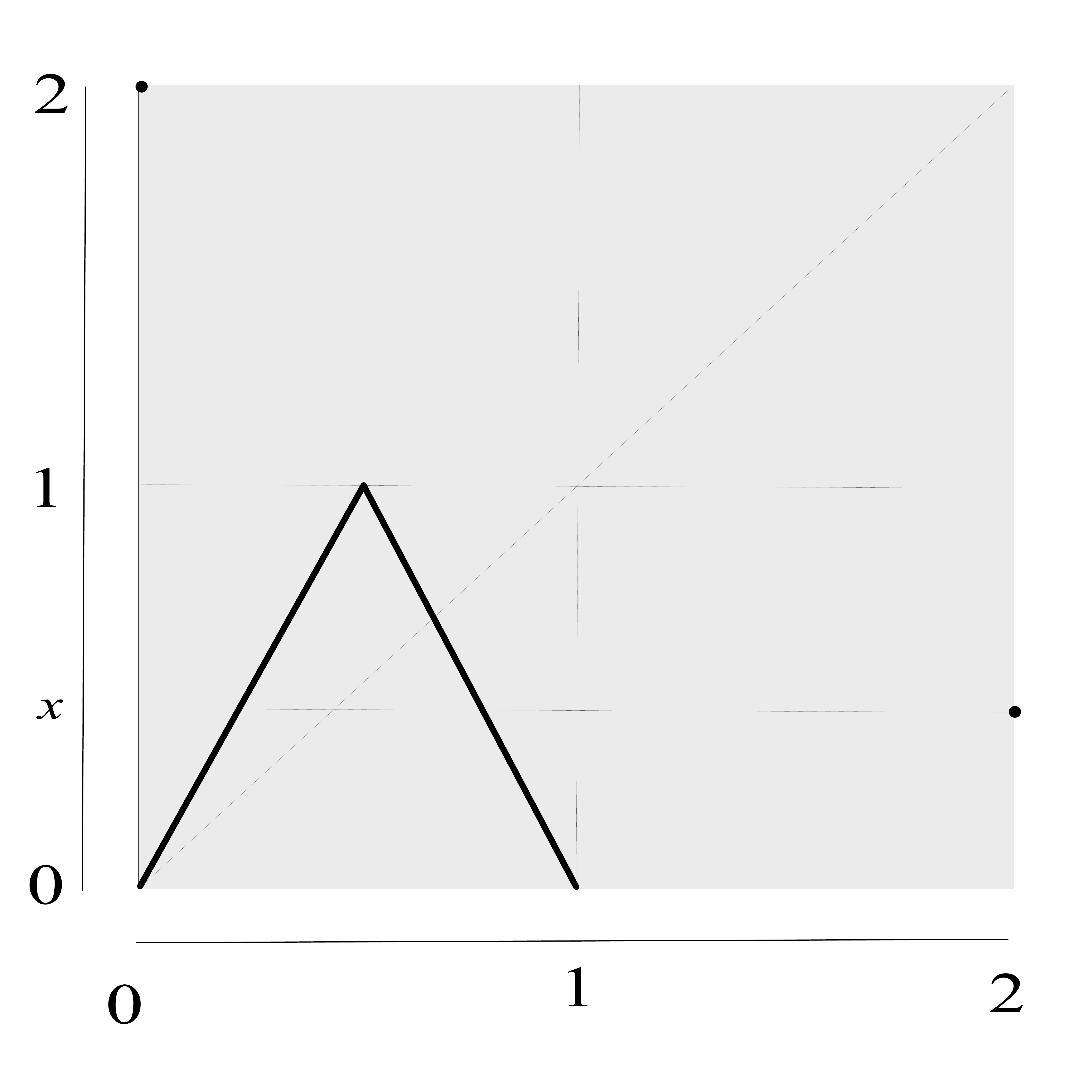}
	\caption{The relation $G$ from Example \ref{exxi}}
	\label{fi}
\end{figure} 
\noindent To see that $(X,G)$ satisfies \ref{133} of Theorem \ref{transitiveCHARACTERIZATION11}, let $U$ be a non-empty open set in $X$.   To see that $\bigcup_{k=0}^{\infty}G^k(U)$ is dense in $X$, let $V$ be a non-empty open set in $X$.  To show that $\bigcup_{k=0}^{\infty}G^k(U)\cap V\neq \emptyset$, we consider the following cases.
\begin{enumerate}
\item $2\in U\cap V$.  Since $U\cap V\neq \emptyset$, it follows that $(\bigcup_{k=0}^{\infty}G^k(U))\cap V\neq \emptyset$.
\item $2\not \in U$ and $2\in V$.  Since $f$ is locally eventually onto, it follows that there is a positive integer $n$ such that $0\in f^n(U)$. Then $0\in G^n(U)$ and it follows that $2\in G^{n+1}(U)$. Therefore,  $(\bigcup_{k=0}^{\infty}G^k(U))\cap V\neq \emptyset$.
\item $2\in U$ and $2\not\in V$.  Then $x\in G(U)$. Since $x\in \tr(f)$, it follows that there is a positive integer $n$ such that $G^n(U)\cap V\neq \emptyset$. Therefore,  $(\bigcup_{k=0}^{\infty}G^k(U))\cap V\neq \emptyset$.
\item $2\not \in U$ and $2\not\in V$.   Since $([0,1],f)$ is transitive, it follows from \cite[Theorem 4.2, page 10]{B} that $(\bigcup_{k=0}^{\infty}f^k(U))\cap V\neq \emptyset$. Therefore,  $(\bigcup_{k=0}^{\infty}G^k(U))\cap V\neq \emptyset$.
\end{enumerate} 
To see that $(X,G)$ does not satisfy \ref{144} of Theorem \ref{transitiveCHARACTERIZATION11}, let $U=\{2\}$. It follows that 
$\Cl\left(\bigcup_{k=1}^{\infty}G^{k}(U)\right)=\Cl\left(\left\{x,f(x),f^2(x),\ldots\right\}\right)=[0,1]\neq X$.
\end{example}
{
\begin{definition}
Let $(X,G)$ be a CR-dynamical system.  We say that $(X,G)$ is $+$transitive, if for all non-empty open sets $U$ and $V$ in $X$,  there is a positive integer $n$ such that 
$$
G^{n}(U)\cap V\neq \emptyset.
$$
\end{definition}
\begin{observation}\label{mojelilike}
Note that the following hold.
\begin{enumerate}
\item For any CR-dynamical system $(X,G)$ it holds that if $(X,G)$ is  $+$transitive, then $(X,G)$ is transitive.
\item There are transitive CR-dynamical systems that are not $+$transitive (see Example \ref{exxi}).
\item Let $(X,G)$ be a CR-dynamical system. It follows from Theorem \ref{transitiveCHARACTERIZATION11}  that if $\isolated(X)=\emptyset$, then
$$
(X,G) \textup{ is transitive } \Longleftrightarrow  (X,G) \textup{ is } \textup{+transitive.} 
$$ 
\item Let $(X,G)$ be a CR-dynamical system. It follows from Theorem \ref{transitiveCHARACTERIZATION11} that 
$$
(X,G) \textup{ is +transitive } \Longleftrightarrow  (X,G^{-1}) \textup{ is } \textup{+transitive.} 
$$ 
\end{enumerate}
\end{observation}
\begin{theorem}
Let $(X,G)$ be a CR-dynamical system. If $(X,G)$ is transitive or $+$transitive, then 
$$
p_1(G)=p_2(G)=X.
$$
\end{theorem}
\begin{proof}
First, suppose that $(X,G)$ is transitive.  If $X$ is degenerate, there is nothing to prove. So, suppose that $X$ is non-degenerate.  First, we prove that $p_1(G)=X$. To show this, let $U=X\setminus p_1(G)$. Suppose that $U\neq \emptyset$.  Then $U$ is a non-empty open subset of $X$. Let $x\in U$ and let $y\in X\setminus \{x\}$. Such a $y$ does exist since $X$ is a non-degenerate space. Let $\varepsilon >0$ be such that 
$$
B(x,\varepsilon)\subseteq U \textup{ and } B(x,\varepsilon)\cap B(y,\varepsilon)=\emptyset.
$$ 
Since $(X,G)$ is transitive, there is a non-negative integer $n$ such that 
$$
G^n(B(x,\varepsilon))\cap B(y,\varepsilon)\neq \emptyset.
$$
Note that since $B(x,\varepsilon)\cap B(y,\varepsilon)=\emptyset$, it follows that $n>0$.  Therefore, $B(x,\varepsilon)\cap p_1(G)\neq \emptyset$, which is a contradiction, since $B(x,\varepsilon) \subseteq U$. It follows that $U=\emptyset$ and, therefore, $p_1(G)=X$.  By Observation \ref{naminusena}, $(X,G^{-1})$ is also transitive.  Using this and what we have just proved, it follows that 
$$
p_2(G)=p_1(G^{-1})=X.
$$
Next, suppose that $(X,G)$ is $+$transitive.  By Observation \ref{mojelilike}, $(X,G)$ is also a transitive CR-dynamical system, therefore $p_1(G)=p_2(G)=X$ follows.
\end{proof}
}
Next, we give an example of a transitive CR-dynamical system  $(X,G)$ such that $X$ has no isolated points and $\trans_2(G)=\emptyset$.  This is an example that proves that the statement of \cite[Theorem 9, page 3]{SS} is incorrect. 

\begin{example}\label{exhura}
Let $X=[0,1]$, and let 
$
f_{1}:\Big[0,\frac{1}{2}\Big]\rightarrow \Big[0,\frac{1}{2}\Big]
$
 be defined by $f_1(t)=2t$ if $t\in[0,\frac{1}{4}]$ and $f_1(t)=1-2t$ if $t\in[\frac{1}{4},\frac{1}{2}]$ and let 
 $
 f_{2}:\Big[\frac{1}{2},1\Big]\rightarrow \Big[\frac{1}{2},1\Big]
 $
  be defined by $f_2(t)=2t-\frac{1}{2}$ if $t\in[\frac{1}{2},\frac{3}{4}]$ and $f_2(t)=\frac{5}{2}-2t$ if $t\in[\frac{3}{4},1]$.   Let $x_1\in \tr(f_1)$ and $x_2\in \tr(f_2)$, and let  $G$ be defined by
$$
G=\Gamma(f_1)\cup \Gamma(f_2)\cup\{(0,x_2),(1,x_1)\},
$$
see Figure \ref{fighura}. Note that $x_1\in (0,\frac{1}{2})$ while $x_2\in (\frac{1}{2},1)$.
\begin{figure}[h!]
	\centering
		\includegraphics[width=20em]{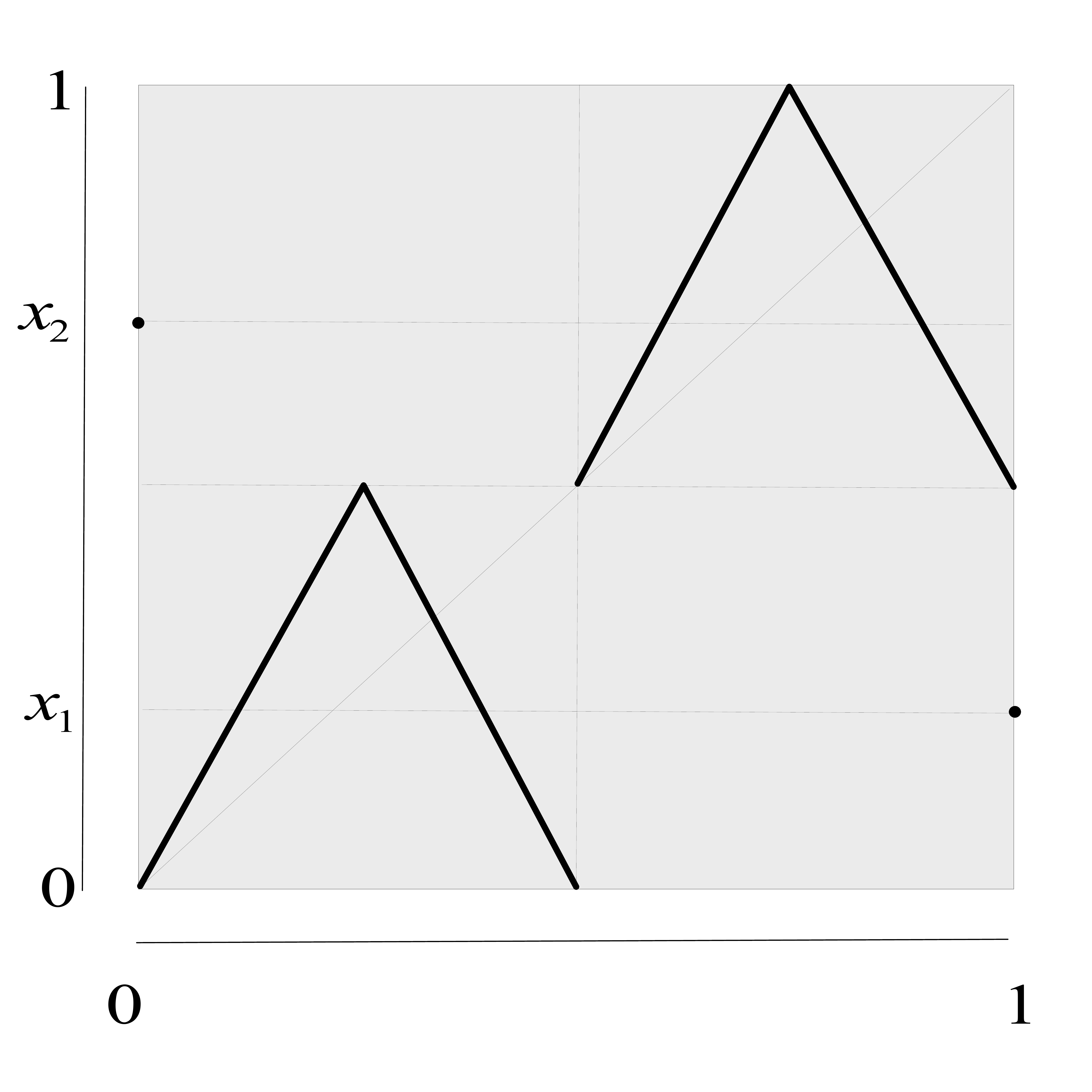}
	\caption{The relation  $G$ from Example \ref{exhura}}
	\label{fighura}
\end{figure} 
 We first show that $(X,G)$ is transitive. Let $U$ and $V$ be non-empty open sets in $X$. 
We consider the following possible cases.
\begin{enumerate}
\item $U\cap [0,\frac{1}{2}]\neq \emptyset$ and $V\cap [0,\frac{1}{2}]\neq \emptyset$. Then there is a non-negative integer $n$ such that $f_1^n(U)\cap V\neq \emptyset$, since $([0,\frac{1}{2}],f_1)$ is transitive.  Let $n$ be such a non-negative integer.  Since $f_1^n(U)\subseteq G^n(U)$ it follows that $G^n(U)\cap V\neq \emptyset$.
\item $U\cap [\frac{1}{2},1]\neq \emptyset$ and $V\cap [\frac{1}{2},1]\neq \emptyset$. Then there is a non-negative integer $n$ such that $f_2^n(U)\cap V\neq \emptyset$, since $([\frac{1}{2},1],f_2)$ is transitive.  Let $n$ be such a non-negative integer.  Since $f_2^n(U)\subseteq G^n(U)$ it follows that $G^n(U)\cap V\neq \emptyset$.
\item  $U\cap [0,\frac{1}{2}]\neq \emptyset$ and $V\cap [\frac{1}{2},1]\neq \emptyset$. Since $f_1$ is locally eventually onto, there is a positive integer $n$ such that $0\in f_1^n(U)$. Let $n$ be such an integer.  Then $x_2\in G^{n+1}(U)$. Since $x_2\in \tr(f_2)$, it follows that there is a non-negative integer $k$ such that $f_2^k(x_2)\in V$.  Let $k$ be such an integer.  It follows that $G^{n+k+1}(U)\cap V\neq \emptyset$.
\item  $U\cap [\frac{1}{2},1]\neq \emptyset$ and $V\cap [0,\frac{1}{2}]\neq \emptyset$. Since $f_2$ is locally eventually onto, there is a positive integer $n$ such that $1\in f_2^n(U)$. Let $n$ be such an integer.  Then $x_1\in G^{n+1}(U)$. Since $x_1\in \tr(f_1)$, it follows that there is a non-negative integer $k$ such that $f_1^k(x_1)\in V$.  Let $k$ be such an integer.  It follows that $G^{n+k+1}(U)\cap V\neq \emptyset$.
\end{enumerate} 
To show that $\trans_2(G)=\emptyset$, suppose that this is not the case and let $x\in \trans_2(G)$.  Let $\mathbf y=(y_1,y_2,y_3,\ldots)\in \star_{i=1}^{\infty}G$ be such that $y_1=x$ and $\orbit_G(\mathbf y)$ is dense in $X$.  We distinguish the following possible cases.
\begin{enumerate}
\item $y_1\leq \frac{1}{2}$.  Let $n$ be the smallest positive integer such that $y_n>\frac{1}{2}$. Then $y_n=x_2$. Since $x_2\in \tr(f_2)$, it follows that for each non-negative integer $k$, $f_2^k(x_2)\neq 1$.  Therefore, for each positive integer $k\geq n$, $y_k\geq \frac{1}{2}$.  It follows that
$$
\Cl(\orbit_G(\mathbf y))=\{y_1,y_2,y_3,\ldots,y_{n-1}\}\cup \left[\frac{1}{2},1\right]=\{0\}\cup \left[\frac{1}{2},1\right]\neq [0,1]
$$
and this is a contradiction since $\orbit_G(\mathbf y)$ is dense in $[0,1]$. 
\item $y_1\geq \frac{1}{2}$.  Let $n$ be the smallest positive integer such that $y_n<\frac{1}{2}$. Then $y_n=x_1$. Since $x_1\in \tr(f_1)$, it follows that for each non-negative integer $k$, $f_1^k(x_1)\neq 0$.  Therefore, for each positive integer $k\geq n$, $y_k\leq \frac{1}{2}$.  It follows that
$$
\Cl(\orbit_G(\mathbf y))=\{y_1,y_2,y_3,\ldots,y_{n-1}\}\cup \left[0,\frac{1}{2}\right]\neq [0,1]
$$
and this is a contradiction since $\orbit_G(\mathbf y)$ is dense in $[0,1]$. 
\end{enumerate}
\end{example}
Next, we give an example of a non-transitive CR-dynamical system  $(X,G)$ such that  $\trans_2(G)\neq\emptyset$.  This is an example that proves that the statement of  \cite[Proposition 8, page 2]{SS} is incorrect.  
\begin{example}
Let $X=\{1,2\}$ and let $G=\{(1,2),(2,2)\}$. Then $\trans_2(G)=\{1\}$.  Let $U=\{2\}$ and $V=\{1\}$. Then $G^n(U)=U$ for each non-negative integer $n$. Therefore,  for each non-negative integer $n$, $G^n(U)\cap V=\emptyset$.  It follows that $(X,G)$ is not transitive. 
\end{example}
We conclude the paper by stating and proving the following theorem.
\begin{theorem}
Let $(X,G)$ be a CR-dynamical system such that $\trans_2(G)\neq \emptyset$.  If $X$ has no isolated points, then  $(X,G)$ is transitive. 
\end{theorem}
\begin{proof}
To see that $(X,G)$ is transitive, let $U$ and $V$ be non-empty open sets in $X$.  Also, let $x\in \trans_2(G)$ and let $\mathbf x=(x_1,x_2,x_3,\ldots)\in \star_{i=1}^{\infty}G$ such that $x_1=x$. Then $\orbit_G(\mathbf x)$ is dense in $X$.  Let $n$ be a positive integer such that $x_n\in U$.  Let $W=V\setminus\{x_1,x_2,x_3,\ldots ,x_n\}$. Then there is a positive integer $m>n$ such that $x_m\in W$. It follows that $G^{m-n}(U)\cap V\neq \emptyset$.
\end{proof}


\noindent I. Bani\v c\\
              (1) Faculty of Natural Sciences and Mathematics, University of Maribor, Koro\v{s}ka 160, SI-2000 Maribor,
   Slovenia; \\(2) Institute of Mathematics, Physics and Mechanics, Jadranska 19, SI-1000 Ljubljana, 
   Slovenia; \\(3) Andrej Maru\v si\v c Institute, University of Primorska, Muzejski trg 2, SI-6000 Koper,
   Slovenia\\
             {iztok.banic@um.si}           
     
				\-
			
	\noindent G.  Erceg\\
             Faculty of Science, University of Split, Rudera Bo\v skovi\' ca 33, Split,  Croatia\\
{{goran.erceg@pmfst.hr}       }    
                 
%
                 \-
			
		\noindent S.  Greenwood\\
            Department of Mathematics, University of Auckland, 38 Princes Street, 1010 Auckland, New Zealand\\
{{s.greenwood@auckland.ac.nz}       }    
%
         
                 	\-
				
  \noindent J.  Kennedy\\
             Department of Mathematics, Lamar University, 200 Lucas Building, P.O. Box 10047, Beaumont, TX 77710 USA\\
{{kennedy9905@gmail.com}       }    
\end{document}